\documentclass{article}
\usepackage{graphicx, amsmath, xcolor, tensor, physics, amsfonts, amssymb, appendix, amsthm}
\usepackage{hyperref}
\hypersetup{
    colorlinks=true,
    linkcolor=blue,
    filecolor=magenta,      
    urlcolor=cyan,
    }
\allowdisplaybreaks

\usepackage{etoolbox}
\makeatletter
\pretocmd\@maketitle
  {\def\@makefnmark{\hbox{\@textsuperscript{\normalfont\@thefnmark}}}}
  {}{\FAIL}
\makeatother

\title{Balanced and Aeppli Parameters for the Heterotic Moduli}
\author{S\'ebastien Picard\thanks{Department of Mathematics, UBC, 1984 Mathematics Road,
    Vancouver, BC, Canada, \href{spicard@math.ubc.ca}{spicard@math.ubc.ca}} and Pei-Lin Wu\thanks{Department of Mathematics, UBC, 1984 Mathematics Road,
    Vancouver, BC, Canada, \href{peilinwu@math.ubc.ca}{peilinwu@math.ubc.ca}}}

\date{\today}

\begin{document}

\newtheorem{theorem}{Theorem}
\newtheorem{lemma}[theorem]{Lemma}

\theoremstyle{remark}
\newtheorem*{rmk}{Remark}

\newcommand{\A}{\mathcal{A}}
\newcommand{\C}{\mathcal{C}}
\newcommand{\D}{\mathcal{D}}
\newcommand{\F}{\mathcal{F}}
\newcommand{\M}{\mathcal{M}}
\renewcommand{\Tr}{\text{Tr}}
\newcommand{\pr}{\text{pr}}
\newcommand{\vol}{\text{vol}}
\newcommand{\id}{\text{id}}
\newcommand{\End}{\text{End}}
\newcommand{\IEnd}{\mathcal{E}nd}
\renewcommand{\Im}{\text{Im}}
\newcommand{\Ker}{\text{ker}}
\newcommand{\rk}{\text{rank}}

\maketitle
\begin{abstract}
   In this paper, we fix the complex structure and explore the moduli space of the heterotic system by considering two different yet ``dual" deformation paths starting from a K\"ahler solution. They correspond to deformation along the Bott-Chern cohomology class and the Aeppli cohomology class respectively. Using the implicit function theorem, we prove the stability of the existence of heterotic solutions under these two deformations and hence establish an initial step in constructing local moduli coordinates around a K\"ahler solution. 
\end{abstract}

\tableofcontents

\section{Introduction}
Calabi-Yau threefolds were proposed as compactifications of extra dimensions in string theory in the work of Candelas-Horowitz-Strominger-Witten \cite{CHSW}. Understanding the landscape of Calabi-Yau threefolds is fundamental for application to string theory, and a study of the local geometry of the parameter space goes back to Candelas and de la Ossa \cite{CdlO}. 
\smallskip
\par The heterotic system is a set of geometric equations derived from heterotic string theory which provides a natural extension of the theory of Calabi-Yau threefolds to non-K\"ahler complex manifolds. It is currently an active field of research to understand its parameter space.
\smallskip
\par Before surveying recent works on this parameter space, we mention that the heterotic system interacts with various fields of pure mathematics. For example, this system can be recast in terms of generalized geometry; see \cite{AGS,OS2014,MGF14} for initial observations and \cite{MGF17,ASTW,MGF-JDG} for follow-up work including a moment map interpretation of the equations. There is also a parabolic version of the system, the Anomaly flow, which can be understood as supersymmetric Ricci flow with higher $\alpha'$ corrections \cite{AMP23,PPZ2,PPZ1,PPZ3,FPPZ22}. Special solutions to the heterotic system are constructed in \cite{DRS,Fei16,FeiYau,FHP,FIUV,FGV,FY,FTY,Grant,OUV}.
\smallskip
\par The framework for understanding variations of families of solutions to the system was initiated by Anderson-Gray-Sharpe \cite{AGS}, Svanes-de la Ossa \cite{OS2014} and Garcia-Fernandez-Rubio-Tipler \cite{MGF17}. There has been much follow-up work in both the string theory and pure mathematics communities. In the string theory literature, we mention e.g. \cite{ADMSS,ASTW,OHS,McOSva}, and in \cite{CdlOMc,CdlOMcS}, the geometry of the total space is investigated and the moduli space is given a K\"ahler structure as a consequence of supersymmetry. These works presuppose the existence of local coordinates on the heterotic moduli and explore the induced geometric structures on the total space.
\smallskip
\par In the mathematics literature on the heterotic moduli space, we mention \cite{MGF17,MGF-JDG} for a natural construction of the infinitesimal moduli by symmetries in generalized geometry, and \cite{MGF20} for a relation to the deformation theory of holomorphic string algebroids.
\smallskip
\par Compared to the works mentioned above, our main focus is the construction of families of solutions. Given a path of parameters in cohomology, we explain how to create associated solutions to the system of equations. We will prove the stability of the existence of heterotic solutions along these two paths, which provide two sets of local moduli coordinates around a given heterotic solution. 
\smallskip
\par Our setup here will assume the existence of a background K\"ahler metric. In this setup, the first solutions to the heterotic system were constructed by Li and Yau \cite{LY05} by the inverse function theorem. This construction was later generalized by Andreas and Garcia-Fernandez \cite{AG2012}. A new proof was given in \cite{CPY2022}, where the inverse function theorem was applied in a fixed balanced class. Compared to \cite{CPY2022}, the present work uses the implicit function theorem to vary the balanced class parameter.
\smallskip
\par Two distinguished cohomology classes can be attached to the heterotic system: the balanced class and the Aeppli class. The roles of these classes are explored in \cite{MGF22,MGF23,MGF-JDG}. We will show how to use these classes as parameters on the moduli, and we construct solutions with prescribed class nearby a K\"ahler class. We state our main results below. 

\subsection*{Outline}
Let $X$ be a complex manifold of dimension 3 with holomorphic volume form $\Omega$ and hermitian metric $\omega$. Let $E \rightarrow X$ be a complex vector bundle with connection $A$. Let $\alpha'>0$. We will construct families of solutions to the following system:
\begin{align} \label{system-of-eqn}
 & \ d(|\Omega|_\omega \omega^2)=0, \quad i \partial \bar{\partial} \omega = \alpha' ({\rm Tr} \, R_{\omega} \wedge R_{\omega} - {\rm Tr} \, F_A \wedge F_A )   \nonumber\\
 & \ F_A^{2,0}=F_A^{0,2}=0, \quad F_A \wedge \omega^2 = 0,
\end{align}
where $F_{A}$ is the curvature of the connection $A$, and $R_{\omega}$ is the curvature of the Chern connection of the metric $\omega$. Furthermore, we make the following assumptions on the background data: [BG]
\begin{align*}
    &\bullet X \text{ admits a K\"ahler Ricci-flat metric } \omega_{\rm CY} \, \cite{Yau78} \tag{BG1} \\
    &\bullet c_1(X) = c_1(E) = 0, \quad c_2(X) = c_2(E) \tag{BG2} \\
    &\bullet E \rightarrow (X, \omega_{\rm CY}) \text{ is a stable holomorphic bundle} \tag{BG3} \\
    &\bullet \text{Def(} E \text{) is smooth and of dimension } h^1(X, \operatorname{End} E) \tag{BG4}
\end{align*}

\begin{rmk}By the Donaldson-Uhlenbeck-Yau theorem \cite{Donaldson,UY}, if $E \rightarrow (X,\omega_{\rm CY})$ is stable then there exists a Hermitian-Yang-Mills metric $h_{\rm DUY}$ such that its Chern connection $A$ solves $F_{A} \wedge \omega_{\rm CY}^2=0$.
  \end{rmk}

  \begin{rmk} We briefly discuss the necessity of the topological constraints (BG2). As $X$ admits a nowhere vanishing holomorphic section $\Omega \in H^0(X,K_X)$, we have $c_1(X)=0$. We require $c_{1}(E)=0$ so that the $\omega$-degree of $E$ is zero. The anomaly cancellation relation in \eqref{system-of-eqn} implies that $c_{2}(X)=c_{2}(E)$.
  \end{rmk}

  \begin{rmk}
In (BG4), Def($E$) denotes the deformation space of holomorphic structures on $E$. The assumption (BG4) is sometimes stated as $E$ having unobstructed deformations \cite{Huybrechts}. Examples of such bundles over complete intersection Calabi-Yau threefolds were studied by Huybrechts \cite{Huybrechts} and include $E=T^{1,0}X$, as well as extensions of $T^{1,0}X$ by $\mathcal{O}_X$ where $E=T^{1,0}X \oplus \mathcal{O}_X$ as smooth bundles.
  \end{rmk}

\par To construct a moduli space of solutions to the equations \eqref{system-of-eqn}, we must first find appropriate moduli parameters. For this, we take the classic work of Candelas-de la Ossa \cite{CdlO} on the local moduli of K\"ahler Ricci-flat metrics as the starting point. There the parameter space is a product of the K\"ahler cone and the space of complex structures. In the present work, we consider a simplified setup where the complex structure of the manifold $X$ is fixed. We are left with the K\"ahler class 
\begin{equation*}
[\omega]\in H^{1,1}(X,\mathbb{R}).
\end{equation*}
By Yau’s theorem \cite{Yau78}, each K\"ahler class $[\omega]$ contains a unique K\"ahler Ricci-flat metric $\omega_{\rm CY} \in [\omega]$.
\smallskip
\par Returning to the system \eqref{system-of-eqn}, it is not $\omega$ that is closed, but rather $|\Omega |_\omega \omega^2$. Accordingly, we write
\begin{equation*}
[| \Omega |_\omega \, \omega^2] \in H^{2,2}(X,\mathbb{R}).
\end{equation*}
On the K\"ahler Calabi-Yau threefold $X$ we have $h^{1,1}=h^{2,2}$, so the dimensions of these parameter spaces are the same. Our goal is to construct exact solutions to \eqref{system-of-eqn} (to all orders in deformation theory) given a parameter $\mathfrak{b} \in H^{2,2}(X,\mathbb{R})$ nearby a K\"ahler class. Since the system \eqref{system-of-eqn} also involves the Hermitian-Yang-Mills equation, we also incorporate a parameter $\alpha \in H^{0,1}({\rm End} \, E)$ which represents the moduli of Hermitian-Yang-Mills connections. Assumption (BG4) associates to each parameter $\alpha$ a holomorphic structure $\bar{\D}_\alpha: \Gamma(E) \rightarrow \Omega^{0,1}(E)$ with $\bar{\D}_\alpha^2=0$ which we also deform to a solution to the system \eqref{system-of-eqn}.
\smallskip
\par We now state our results more precisely, and refer to Bott-Chern cohomology $H^{2,2}_{\rm BC}$ as this may be the correct notion of cohomology for parametrizing solutions on non-K\"ahler manifolds; in the present case $X$ admits a background K\"ahler metric and hence $H^{2,2}_{\rm BC}$ is equal to the usual Dolbeault cohomology $H^{2,2}$.
\smallskip
\par Let $\mathcal{M}$ be the space of all hermitian metrics on $X$. Let $\mathcal{J}_E$ be the space of all connections on $E$. For $\epsilon>0$, we construct an injective map
\begin{equation*}
    \sigma: B_\epsilon(0) \subset H^{2,2}_{\rm BC}(X,\mathbb{R}) \times H^{0,1}({\rm End} \, E) \rightarrow \mathcal{M} \times \mathcal{J}_E
\end{equation*}
denoted
\begin{equation*}
    \sigma(\mathfrak{b},\alpha) = (\omega_{\mathfrak{b},\alpha}, A_{\mathfrak{b},\alpha})
\end{equation*}
such that $(\omega_{\mathfrak{b},\alpha}, A_{\mathfrak{b},\alpha})$ solves the system \eqref{system-of-eqn} and
\begin{equation*}
  [| \Omega |_{\omega_{\mathfrak{b},\alpha}} \, \omega_{\mathfrak{b},\alpha}^2] = [|\Omega|_{\omega_{\rm CY}} \omega^2_{\rm CY}] + \mathfrak{b} \in H^{2,2}_{\rm BC}(X,\mathbb{R}).
\end{equation*}
Setting the balanced class parameter $\mathfrak{b}=0$, we obtain $H^{0,1}({\rm End} \, E)$ solutions to the system with balanced class given by the reference Calabi-Yau metric; this generalizes previous work of \cite{CPY2022}. In summary, we construct a family of solutions locally parametrized by
\begin{equation*}
     H^{2,2}_{\rm BC}(X,\mathbb{R}) \times H^{0,1}({\rm End} \, E).
\end{equation*}
Returning to the analogy with K\"ahler geometry, the first piece is similar to the K\"ahler class moduli $H^{2,2}(X,\mathbb{R}) \cong H^{1,1}(X,\mathbb{R})$ which parametrize K\"ahler Ricci-flat metrics. The second piece is the moduli of Hermitian-Yang-Mills connections $H^{0,1}({\rm End} \, E)$ with fixed K\"ahler metric.

\begin{theorem}[Stability of existence under Bott-Chern deformation]\label{thm-bc}
Suppose the background assumptions [BG] are satisfied. There exists a small parameter $\alpha'>0$ such that the heterotic system \eqref{system-of-eqn} admits solutions $(\Tilde{\omega},\Tilde{A})$ nearby a K\"ahler structure $(\omega_{\rm CY}, A)$ along deformation paths parameterized by 
\begin{equation*}
[\mathfrak{b}]\in H^{2,2}_{\rm BC}(X,\mathbb{R}), \quad [\alpha_1]\in H^{0,1}(X, {\rm End} \, E).    
\end{equation*}
This equips the heterotic moduli with local coordinates $([\mathfrak{b}],[\alpha_1])$ in a small $\epsilon$-neighborhood of a K\"ahler solution.
\end{theorem}

There is another cohomology class associated to a solution to the system \eqref{system-of-eqn}. This alternate analog to the K\"ahler class was introduced in \cite{MGF22} and is an element of $H^{1,1}_{\rm A}(X,\mathbb{R})$ in Aeppli cohomology. There is a slight technical complication where to state the results in this setup, we introduce an additional variable $\theta$ into the system. In this setup, we let $\theta$ be a Hermitian-Yang-Mills connection on the smooth tangent bundle $T^{1,0}X$ and use it to compute ${\rm Tr} \, R_\theta \wedge R_\theta$. By work of Huybrechts \cite{Huybrechts}, complete intersection Calabi-Yau threefolds have the property that $T^{1,0}X$ has unobstructed deformations, hence our assumption (BG4) applies. This alternate setup is commonly used in the literature; see e.g. \cite{AG2012,OS2014,G2016,MGF23}. The Aeppli class of a solution is then denoted
\begin{equation*}
    \underline{a}(\omega,A,\theta) \in H^{1,1}_A(X,\mathbb{R})
\end{equation*}
and we recall the definition in Section \ref{Aeppli-defn}.
\smallskip
\par We can now compare Theorem \ref{thm-bc} with the deformation theory of Aeppli classes developed in \cite{MGF22}. In this alternate setup, there is an injective map
\begin{equation*}
    \sigma: B_\epsilon(0) \subset H^{1,1}_{\rm A}(X,\mathbb{R}) \times H^{0,1}({\rm End} \, E) \times H^{0,1}({\rm End} \, T^{1,0}X) \rightarrow \mathcal{M} \times \mathcal{J}_E \times \mathcal{J}_{T^{1,0}X}
\end{equation*}
denoted
\begin{equation*}
    \sigma(\mathfrak{a},\alpha_1,\alpha_2) = (\omega_{\mathfrak{a},\alpha}, A_{\mathfrak{a},\alpha}, \theta_{\mathfrak{a},\alpha})
\end{equation*}
such that $(\omega_{\mathfrak{b},\alpha}, A_{\mathfrak{b},\alpha})$ solves the system \eqref{system-of-eqn} and
\begin{equation*}
    \underline{a}(\omega_{\mathfrak{a},\alpha}, A_{\mathfrak{a},\alpha}, \theta_{\mathfrak{a},\alpha})  = [\omega_{\rm CY}] + \mathfrak{a}.
\end{equation*}
In other words, the Aeppli class can also be used as a coordinate on the moduli space. At $\mathfrak{a}=0$, there is a $H^{0,1}({\rm End} \, E) \times H^{0,1}({\rm End} \, T^{1,0}X)$ family of solutions all inside the same Aeppli class as the given reference Calabi-Yau metric $\omega_{\rm CY}$. Nearby $[\omega_{\rm CY}]$, there are solutions in a given Aeppli class. 

\begin{theorem}[Stability of existence under Aeppli deformation]\label{thm-a}
   Suppose the background assumptions $[BG]$ are satisfied and $T^{1,0}X$ has unobstructed deformations. There exists a small parameter $\alpha'>0$ such that the heterotic system \eqref{system-of-eqn} with an additional spurious gauge field admits solutions $(\Tilde{\omega},\Tilde{A},\Tilde{\theta})$ nearby a K\"ahler structure $(\omega_{\rm CY},A, \Gamma_{\rm CY})$ along deformation path parameterized by 
   \begin{equation*}
   [\mathfrak{a}]\in H^{1,1}_{\rm A}(X,\mathbb{R}), \quad [\alpha_1]\in H^{0,1}(X, {\rm End} \, E), \quad [\alpha_2]\in H^{0,1}(X, {\rm End} \, T^{1,0}X).
   \end{equation*}
   Thus the heterotic moduli admits the local coordinates $([\mathfrak{a}],[\alpha_1],[\alpha_2])$ in a small $\epsilon$-neighborhood of a K\"ahler solution.
\end{theorem}

Theorem \ref{thm-a} was previously known in the literature, and follows from Corollary 5.14 in \cite{MGF22} by Garcia-Fernandez, Rubio, Shahbazi and Tipler, together with the main result of Andreas and Garcia-Fernandez \cite{AG2012}. We present here a self-contained exposition which is streamlined such that both the balanced and Aeppli deformations can be compared in parallel.

\subsection*{Acknowledgements}
We thank M. Garcia-Fernandez, J. McOrist, and E. Svanes for helpful discussions and useful comments. We also thank the referee for a careful reading and helpful suggestions.

\subsection{Notation and conventions}
\subsubsection{The base manifold $X$}
Let $X$ be a compact complex manifold of complex dimension $n=3$ which admits a nowhere vanishing holomorphic 3-form $\Omega\in\Lambda^{3,0}(X,\mathbb{C})$. Let $\omega\in\Lambda^{1,1}(X,\mathbb{R})$ be a hermitian metric on $X$. 

If $\{x^{\alpha}\}_{1}^{3}$ are the holomorphic coordinates of $X$, then $\omega$ can be locally expressed by
\begin{equation*}
\omega=ig_{\alpha\bar{\beta}}dx^{\alpha}\wedge d\bar{x}^{\beta}=ig_{\alpha\bar{\beta}}dx^{\alpha\bar{\beta}},
\end{equation*}
where $g_{\alpha\bar{\beta}}$ is hermitian ($\overline{g_{\alpha\bar{\beta}}}=g_{\beta\bar{\alpha}}$) and positive-definite. The inverse of $g_{\alpha\bar{\beta}}$ is denoted by $g^{\bar{\alpha}\beta}$, and $g_{\alpha\bar{\gamma}}g^{\bar{\gamma}\beta}=\tensor{\delta}{_{\alpha}^{\beta}}$. 

The nowhere vanishing holomorphic 3-form $\Omega$ can be locally expressed by
\begin{equation*}
    \Omega=fdx^1\wedge dx^2\wedge dx^3,
\end{equation*}
where $f$ is a non-vanishing holomorphic function. The norm of $\Omega$ with respect to $\omega$, denoted as $|\Omega|_{\omega}$ is defined via the relation $i\Omega\wedge\bar{\Omega}=|\Omega|_{\omega}^2\tfrac{\omega^3}{3!}$. It has a local expression
\begin{equation}\label{eq:localOmega}
    |\Omega|_{\omega}^2=f\bar{f}(\det g)^{-1}.
\end{equation}

From a hermitian metric $\omega$, we can form the Chern connection $\Gamma$, locally given by $\Gamma^{\kappa}_{\alpha\beta}=g^{\kappa\bar{\gamma}}\partial_{\alpha}g_{\beta \bar{\gamma}}$. The corresponding Chern curvature $R$ is locally given by $\tensor{R}{_{\alpha\bar{\beta}}^{\mu}_{\nu}}=-\partial_{\bar{\beta}}\Gamma^{\mu}_{\alpha\nu}$.

The hermitian metric $\omega$ is called K\"ahler if $d\omega=0$, and called conformally balanced if 
\begin{equation*}
d(|\Omega|_{\omega}\omega^2)=0.
\end{equation*}
If $X$ admits a K\"ahler metric $\omega$, then by Yau's theorem \cite{Yau78} there exists a unique K\"ahler Ricci-flat metric $\omega_{\mathrm CY} \in [\omega]$. K\"ahler Ricci-flat metrics are conformally balanced as in this case $|\Omega|_\omega$ is constant.

The Lefschetz operator $L_{\omega}$ is defined by $L_{\omega}\gamma=\omega\wedge\gamma$ for any form $\gamma$. The contraction operator $\Lambda_{\omega}=\star_{\omega}L_{\omega}\star_{\omega}$ can locally expressed by $i\Lambda_{\omega}\gamma=g^{\bar{\beta}\alpha}\gamma_{\alpha\bar{\beta}}$ for any $2$-form $\gamma$, such that $\Lambda_{\omega}\omega=n=3$. The counting operator $H$ is defined to be $H=\sum_{k=0}^{2n}(k-n)\Pi^{k}$, with $\Pi^{k}$ the natural projection to $k$-forms. Then, the action of $L_{\omega}$, $\Lambda_{\omega}$, and $H$ defines natural $\mathfrak{sl}_{2}$-representations, in particular, we have the following commutation relations (see Proposition 1.2.26 from \cite{D2015}),
\begin{equation}\label{eq:KahlerId}
    [H, L_{\omega}]=2L_{\omega},\quad [H,\Lambda_{\omega}]=-2\Lambda_{\omega}, \quad[L_{\omega},\Lambda_{\omega}]=H.
\end{equation}

\subsubsection{The gauge bundle $E$ over $X$}
Let $E\to X$ be a smooth complex vector bundle over $X$ of rank $r$ equipped with a hermitian metric $h$. Our convention for the inner product of sections $x,y \in \Gamma(E)$ is 
\[
\langle x,y \rangle_{h} = y^* h x, \quad \langle x,y \rangle_{h} = \overline{y^\alpha} h_{\bar{\alpha} \beta} x^\beta,
\]
where in a local frame $\{ e_\alpha \}$, we expand $x = x^\alpha e_\alpha$ and $h = h_{\bar{\alpha} \beta} \, e^\beta \otimes \overline{e^\alpha}$. 

Let $D=d+A$ be the covariant derivative of a connection $A$ on $E$ which acts on sections $s \in \Gamma(E)$ by
\[
D s = d s + A s,
\]
and acts on endomorphisms $\gamma \in \Gamma({\rm End} \, E)$ by
  \[
D^{{\rm End} E} \gamma = d \gamma + [A,\gamma].
    \]
    For simplicity, we will from now on simply write the induced connection on ${\rm End} \, E$ acting on $\gamma \in \Gamma({\rm End} \, E)$ by $D \gamma$ rather than $D^{{\rm End} E} \gamma$.

    The curvature of the connection $A$ is denoted by $F_{A}=dA+A\wedge A$. We may decompose the connection in the following way
\[
D = \D + \bar{\D}, \quad \D: \Gamma(E) \rightarrow \Omega^{1,0}(E), \quad \bar{\D}: \Gamma(E) \rightarrow \Omega^{0,1}(E).
\]
Suppose $\bar{\D}^2=F_A^{0,2}=0$. In this case, by a well-known integrability result (see e.g. \cite{DonaldsonKronheimer}), $E$ can be given trivializations by holomorphic frames so that $E \rightarrow X$ is a holomorphic vector bundle. In a holomorphic frame, then $\bar{\D} = \bar{\partial}$. If we require $D$ to be metric compatible with $h$, meaning
\begin{equation} \label{metric-comp}
\partial \langle x,y \rangle_h = \langle \D x, y \rangle_h + \langle x, \bar{\D} y \rangle_h,
\end{equation}
then we find in a holomorphic frame that $A = h^{-1} \partial h$. The curvature is then $F_A= F_h=\bar{\partial} (h^{-1} \partial h)$.

Given a holomorphic bundle $E \rightarrow (X,\omega)$ over a K\"ahler manifold, we define the $\omega$-degree of $E$ by 
\begin{equation*}
    \deg(E)=\int_{X}c_1(E)\wedge \omega^{2},
\end{equation*}
and the $\omega$-slope of $E$ by $\mu(E)=\deg(E)/\rk(E)$. Then $E$ is called stable if for any coherent subsheaf $\F \subset \mathcal{O}(E)$ satisfying $0<\rk(F)<\rk(E)$, we have $\mu(\F)<\mu(E)$. The bundle $E$ is said to be poly-stable if $E$ is the direct sum of stable vector bundles with the same $\omega$-slope. The Donaldson-Uhlenbeck-Yau theorem \cite{Donaldson,UY} states that for a compact K\"ahler manifold $(X,\omega)$ and holomorphic vector bundle $E\to X$, the bundle $E$ is poly-stable if and only if $E$ admits a Hermitian-Yang-Mills metric $h$, meaning that its curvature $F_h=\bar{\partial}(h^{-1}\partial h)$ solves $\Lambda_{\omega}F_{h}=\mu\mathbf{1}_{E}$. 

In our setup, we will take as initial data a K\"ahler Calabi-Yau threefold $X$ and a holomorphic bundle $E \to X$ which is stable with respect to $\omega_{\rm CY}$ and satisfies $c_1(E)=0$ and $c_2(E)=c_2(X)$. It is well-known that since $E$ is stable, then it does not admit any holomorphic endomorphism other than multiplies of the identity (see e.g. \cite{SiuBook}) and we will use this later on to rule out elements of the kernel of a linearized operator. We apply the Donaldson-Uhlenbeck-Yau theorem and take the bundle metric $h$ to be Hermitian-Yang-Mills, so that
\begin{equation} \label{ref-metric}
    F_h \wedge \omega_{\rm CY}^2 = 0.
\end{equation}
The Chern connection of $h$ is $D=d+A$ with $A=h^{-1} \partial h$ in a holomorphic frame. The assumption $c_1(E)=0$ gives
\begin{equation*}
    {\rm Tr} \, iF_A = d \gamma
\end{equation*}
for some $\gamma \in \Lambda^1(X)$.

\subsubsection{The Strominger system}
Let $X$ be a complex threefold, $\Omega$ a holomorphic volume form, and $(E,h) \rightarrow X$ a smooth complex vector bundle with metric $h$. With the data $(X,\Omega)$, $(E,h)$ considered fixed, the Strominger system \cite{Strominger} can be translated into a system for a pair $(\omega,A)$ satisfying the following three conditions: 
\begin{enumerate}
    \item The hermitian metric $\omega$ is conformally balanced: $d(|\Omega|_{\omega}\omega^2)=0$,
    \item The unitary connection $A$ solves the Hermitian-Yang-Mills (HYM) equation: $F_{A}^{0,2}=0$ and $F_{A}\wedge \omega^2=0$,
    \item The anomaly cancellation relation holds: 
    \begin{equation} \label{anomaly-eqn}
    i\partial\bar{\partial}\omega=\alpha'(\Tr \, R_{\omega}\wedge R_{\omega} - \Tr \, F_{A}\wedge F_{A}).
    \end{equation}
 
\end{enumerate}

These equations of heterotic systems have been historically split into two terms:
\begin{enumerate}
    \item D-terms \cite{OHS}: $d(|\Omega|_{\omega}\omega^2)=0$ and $F_{A}\wedge \omega^2=0$,
    \item F-terms \cite{AGLO}: $i\partial\bar{\partial}\omega=\alpha'(\Tr R_{\omega}\wedge R_{\omega} - \Tr F_{A}\wedge F_{A})$ and $F_{A}^{0,2}=0$.
\end{enumerate}
Each relates to different supersymmetric couplings in the 4-dimensional effective supergravity theory. 

It is worth noting that throughout the paper, we choose the Chern connection and its corresponding Chern curvature $R_{\omega}$. In the physics literature, one usually takes the Hull connection \cite{H1986} (see \eqref{hull-connection} and \eqref{bismut-connection} in the appendix for the definition) involving the 3-form $H$ and views \eqref{anomaly-eqn} not as an equality but as an expansion where the remaining terms are of order $O(\alpha'^2)$. Nevertheless, since our setup involves a perturbation of a K\"ahler background, our constructed solutions solve the equations of heterotic string theory to appropriate order in $\alpha'$. See Appendix \ref{app:justforchern} for justification for taking the Chern connection and a discussion of the physical system from string theory.

To understand the moduli space of the heterotic system, we can start with a K\"ahler Calabi-Yau solution $d\omega_{\rm CY}=0$ that satisfies the heterotic system at $\alpha'=0$, and probe nearby local solutions. To do this, one has the freedom to fix the deformation ansatz to satisfy some of the three conditions and solve for the remaining, and different partition of the three conditions leads to different deformation ansatzes. The two most natural ansatzes are as follows:
\begin{enumerate}
    \item Bott-Chern case: fix the deformation ansatz to satisfy $d(|\Omega|_{\omega}\omega^2)=0$ and $F^{0,2}_{A}=0$, and solve for $i\partial\bar{\partial}\omega=\alpha'(\Tr R_{\omega}\wedge R_{\omega} - \Tr F_{A}\wedge F_{A})$ and $F_{A}\wedge \omega^2=0$. This approach corresponds to deforming the metric along the Bott-Chern cohomology class. 
    \item Aeppli case: fix the deformation ansatz to satisfy $i\partial\bar{\partial}\omega=\alpha'(\Tr R_{\omega}\wedge R_{\omega} - \Tr F_{A}\wedge F_{A})$ and $F^{0,2}_{A}=0$, and solve for $d(|\Omega|_{\omega}\omega^2)=0$ and $F_{A}\wedge \omega^2=0$. This approach corresponds to deforming the metric along the Aeppli cohomology class. 
\end{enumerate}
It is worth noting that by the isomorphism between the de-Rham cohomology, Bott-Chern cohomology, and Aeppli cohomology on K\"ahler background \cite{DGMS1975}, one can regard these two deformation ansatzes as ``dual" to each other on the initial K\"ahler solution. We hence provide a parallel treatment of these two approaches towards constructing local coordinates on the heterotic moduli space.

\section{Deformations}
Our setup is as follows. We start from a K\"ahler Ricci-flat metric $\omega$ and a Hermitian-Yang-Mills connection $A$, so that $(\omega,A)$ solves the heterotic system at $\alpha'=0$. 
\begin{equation*}
d(|\Omega|_\omega \omega^2)=0, \quad F_A^{0,2}=0, \quad F_A \wedge \omega^2 =0, \quad i \partial \bar{\partial} \omega = 0.    
\end{equation*}
We would like to deform this structure to a nearby solution $(\Tilde{\omega},\Tilde{A})$ solving the system to next order in $\alpha'$, so that the anomaly equation becomes
\begin{equation} \label{anomaly-eqn2}
i\partial\bar{\partial} \Tilde{\omega}=\alpha'(\Tr \, R \wedge R - \Tr \, F_{\tilde{A}}\wedge F_{\tilde{A}}).
\end{equation}
We will deform in two directions. For the first, we deform along the balanced class so that 
\begin{equation*}
    \underline{b}(\tilde{\omega},\tilde{A}) = [|\Omega|_{\omega} \omega^2] +[\mathfrak{b}]
\end{equation*}
where $\mathfrak{b} \in H^{2,2}_{\text{BC}}(X,\mathbb{R})$ is a Bott-Chern class parameter. For the second approach, we deform along the Aeppli class so that
\begin{equation*}
    \underline{a}(\tilde{\omega},\tilde{A}) = [\omega] +[\mathfrak{a}]
\end{equation*}
where $\mathfrak{a} \in H^{1,1}_{\text{A}}(X,\mathbb{R})$ is an Aeppli class parameter. Details and definitions are given below. We note that there is an extra complication in the Aeppli case where we introduce an extra variable $\theta$ into the problem which is used to compute ${\rm Tr} \, R_\theta \wedge R_\theta$. 

These two deformations are non-K\"ahler analogs of varying the K\"ahler class of a Calabi-Yau metric, since on a K\"ahler manifold there holds $H^{2,2}_{\text{BC}}(X,\mathbb{C}) \cong H^{2,2}_{\bar{\partial}}(X,\mathbb{C})$ and $H^{1,1}_{\text{A}}(X,\mathbb{C}) \cong H^{1,1}_{\bar{\partial}}(X,\mathbb{C})$, and the parameter space of the K\"ahler moduli has dimension $h^{1,1} = \dim H^{1,1}_{\bar{\partial}} = \dim H^{2,2}_{\bar{\partial}}$. On a non-K\"ahler manifold, there is nonetheless a duality between $H^{2,2}_{\rm BC}(X,\mathbb{C})$ and $H^{1,1}_{\rm A}(X,\mathbb{C}) $ in analogy to the Poincar\'e duality in the de Rham cohomology, seen by the following pairing
\begin{equation*}
    H^{2,2}_{\rm BC}(X,\mathbb{C})\times H^{1,1}_{\rm A}(X,\mathbb{C})\to \mathbb{C},\quad ([\mathfrak{b}]_{\rm BC},[\mathfrak{a}]_{\rm A})\mapsto\int_{X}\mathfrak{b}\wedge\mathfrak{a},
\end{equation*}
which is known to be well-defined and non-degenerate. 


\subsection{Deformation of metric}
There are two parallel deformation ansatzes corresponding to our two choices of partitioning the equations of the heterotic system. 

\subsubsection{Bott-Chern case} \label{Bott-Chern-defn}
Since we would like our ansatz to satisfy the balanced condition $d(|\Omega|_{\omega}\omega^2)=0$, we therefore choose the deform the metric along the Bott-Chern cohomology class $H^{2,2}_{\text{BC}}(X,\mathbb{C})$. Recall
\begin{equation*}
H^{2,2}_{\text{BC}}(X,\mathbb{C})= \frac{(\ker d) \cap \Lambda^{2,2}}{ ({\rm Im} \, \partial \bar{\partial}) \cap \Lambda^{1,1}}.
\end{equation*}
From a conformally balanced metric $\tilde{\omega}$, we can produce a balanced class
\begin{equation*}
\underline{b}(\tilde{\omega}) = [|\Omega|_{\tilde{\omega}} \, \tilde{\omega}^2] \in  H^{2,2}_{\text{BC}}(X,\mathbb{R}).   
\end{equation*}
With $\omega=\omega_{\rm CY}$ our initial reference metric, we take the deformation for the metric to be of the form
 \begin{equation*}
\underline{b}(\tilde{\omega}) = [|\Omega|_{\omega_{\rm CY}} \omega_{\rm CY}^2] +[\mathfrak{b}].
    \end{equation*}
More precisely, our deformation ansatz is
    \begin{equation}\label{eq:BCPath}
        |\Omega|_{\Tilde{\omega}} \, \Tilde{\omega}^2= |\Omega|_\omega \omega^2+\mathfrak{b}+\Theta,
    \end{equation}
    where 
    \begin{equation*}
(\mathfrak{b}+\Theta)\in\Lambda^{2,2}(X,\mathbb{R}):d\mathfrak{b}=0, \ \Theta=i\partial\bar{\partial}\gamma\in\partial\bar{\partial}\Lambda^{1,1}(X,\mathbb{R}).
    \end{equation*}
For $(\mathfrak{b}+\Theta)$ small enough, taking the square root defines a positive hermitian metric $\tilde{\omega}$; see \cite{Michelsohn, PPZ2} for the precise formula. Then one can see that the deformed metric satisfies the balanced condition,
    \begin{equation*}
        d(|\Omega|_{\Tilde{\omega}} \Tilde{\omega}^2)=d(|\Omega|_{\omega}\omega^2)+d\mathfrak{b}+d(i\partial\bar{\partial}\gamma)=0.
    \end{equation*}
To solve the system, we will substitute this ansatz into the remaining equations.
    
\subsubsection{Aeppli case} \label{Aeppli-defn}
In this setup,  we demand our ansatz to satisfy the anomaly cancellation relation and choose to deform the metric along the Aeppli cohomology class $H^{1,1}_{\text{A}}(X,\mathbb{C})$. Recall
\begin{equation*}
H^{1,1}_{\text{A}}(X,\mathbb{C})= \frac{(\ker \partial \bar{\partial}) \cap \Lambda^{1,1}}{ ({\rm Im} \, \partial) \cap \Lambda^{0,1} \oplus ({\rm Im} \, \bar{\partial}) \cap \Lambda^{1,0} }.
\end{equation*}
    There is one caveat: due to the complexity of $R_{\omega}$, we introduce a spurious degree of freedom, namely an extra gauge field $\theta$ on $T^{1,0}(X)$, to compute $\Tr R\wedge R$ with a modified anomaly cancellation relation,
    \begin{equation} \label{anomaly-eqn3}
     i\partial\bar{\partial}\Tilde{\omega}=\alpha'(\Tr \, \Tilde{R}_{\Tilde{\theta}}\wedge \Tilde{R}_{\Tilde{\theta}} - \Tr \, \Tilde{F}_{\Tilde{A}}\wedge \Tilde{F}_{\Tilde{A}}),
    \end{equation}
    and requiring this spurious gauge field $\tilde{\theta}$ to satisfy the Hermitian-Yang-Mills equation,
    \begin{equation*}
        R_{\tilde{\theta}}\wedge \omega^2=0, \quad R^{0,2}_{\tilde{\theta}}=0.
    \end{equation*}
    Such an approach with the introduction of the gauge field $\theta$ was considered by previous works \cite{AG2012,OS2014,G2016} and others. We discuss why this setup in the K\"ahler large radius limit also solves the physical equations in Appendix \ref{app:justforchern}.

    The topological condition on Chern classes $c_1(E)=c_1(X)=0$ and $c_2(E)=c_{2}(X)$, together with the $\partial \bar{\partial}$-lemma, implies the existence of $\beta\in\Lambda^{1,1}(X,\mathbb{R})$ such that,
    \begin{equation*}
\Tr \, R_{\theta}\wedge R_{\theta} -        \Tr \, F_{A}\wedge F_{A}=i\partial\bar{\partial}\beta.
    \end{equation*}
Here $(\omega,A,\theta)$ are background reference fields. To ensure $\beta$ is real, we may take the real part if necessary since $i \partial \bar{\partial} {\rm Im} \, \beta = 0$.
\smallskip
\par As noticed in \cite{MGF22}, a solution $(\tilde{\omega},\tilde{A},\tilde{\theta})$ to the anomaly cancellation relation \eqref{anomaly-eqn3} creates an Aeppli class \footnote{For this class to be well-defined independently of the choice $(A,\theta,\beta)$, the form $\beta$ must be properly defined; we refer to \cite{PicSurvey} for further details.}
\begin{equation*}
    \underline{a}(\tilde{\omega},\tilde{A},\tilde{\theta}) = \bigg[ \tilde{\omega} - \alpha' R_2[\tilde{\theta},\theta] + \alpha' R_2[\tilde{A},A] - \alpha' \beta \bigg] \in H^{1,1}_{\text{A}}(X,\mathbb{R}).
\end{equation*}
    Here the $R_{2}$'s are the Bott-Chern-Simons secondary characteristics defined respectively by,
    \begin{align}
        i\partial\bar{\partial}R_{2}[\Tilde{A},A]&=\Tr \, \Tilde{F}_{\Tilde{A}}\wedge \Tilde{F}_{\Tilde{A}}- \Tr \, F_{A}\wedge F_{A}, \nonumber\\
        i\partial\bar{\partial}R_{2}[\Tilde{\theta},\theta]&=\Tr \, \Tilde{R}_{\Tilde{\theta}}\wedge \Tilde{R}_{\Tilde{\theta}}- \Tr \, R_{\theta}\wedge R_{\theta}. \nonumber
    \end{align}
    We refer to e.g. \cite{Donaldson} for a construction of the $R_2$. Our construction of $R_2[\tilde{A},A]$ will be slightly different since for $\underline{a}$ to be independent of the choice of references $(A,\theta,\beta)$, one must define $R_2$ as done below in \eqref{defn-R2}; see \cite{PicSurvey} for details.
    \smallskip
    \par Recall that we start from a reference $(\omega,A,\theta)$ where $\omega= \omega_{\rm CY}$ is K\"ahler Ricci-flat, $A$ is Hermitian-Yang-Mills, and $\theta$ is the Chern connection of $\omega$. This means that $(\omega,A,\theta)$ satisfies the heterotic equations  with $\alpha'=0$. We will deform the Aeppli class by
    \begin{equation*}
        \underline{a}(\tilde{\omega},\tilde{A},\tilde{\theta}) = [\omega_{\rm CY}] + [\mathfrak{a}].
    \end{equation*}
    More precisely, we choose the deformation ansatz of the metric to be,
    \begin{equation}\label{eq:APath}
    \Tilde{\omega}=\omega+\mathfrak{a}+b+\alpha'\qty( -R_{2}[\Tilde{A},A]+R_{2}[\Tilde{\theta},\theta] + \beta),
    \end{equation}
    where 
    \begin{equation*}
        (\mathfrak{a}+b)\in \Lambda^{1,1}(X,\mathbb{R}):\partial\bar{\partial}\mathfrak{a}=\partial \bar{\partial} b= 0, \ \ b\in(\Im \, \partial+\Im \, \bar{\partial})\cap\Lambda^{1,1}(X,\mathbb{R}).
\end{equation*}
    Then this ansatz automatically implements the anomaly cancellation relation of the deformed metric \eqref{anomaly-eqn3}. To solve the system, we must substitute this ansatz into the conformally balanced equation. We note that the ansatz \eqref{eq:APath} was introduced in \cite{MGF22} in relation to the notion of positive metric on a Bott-Chern algebroid.
\smallskip
\par We now define $R_{2}[\Tilde{A},A]$ in ansatz \eqref{eq:APath}. First, by Chern-Weil theory, there exists $\theta$ such that
\begin{equation*}
    \Tr \, \Tilde{F}_{\Tilde{A}}\wedge \Tilde{F}_{\Tilde{A}}- \Tr \, F_{A}\wedge F_{A} = d \theta.
\end{equation*}
Applying the $\partial \bar{\partial}$-lemma, we may write $d \theta = i \partial \bar{\partial} \chi$. The issue with $\chi$ is that its construction is ad hoc and its dependence on the unknown connection $\tilde{A}$ is not clear. To give a canonical construction, we use the Kodaira-Spencer \cite{KodairaSpencer} operator
\begin{equation*}
    E = \partial \bar{\partial} \bar{\partial}^\dagger \partial^\dagger + \bar{\partial}^\dagger \partial^\dagger \partial \bar{\partial} + \bar{\partial}^\dagger \partial \partial^\dagger \bar{\partial} +  \partial^\dagger \bar{\partial} \bar{\partial}^\dagger \partial + \bar{\partial}^\dagger \bar{\partial} + \partial^\dagger \partial
\end{equation*}
with adjoints with respect to $g_{\rm CY}$. The operator $E$ is a 4th order self-adjoint elliptic operator, and by elliptic theory, there exists a unique solution $\gamma \in (\ker E)^\perp$ solving
\begin{equation} \label{kodaira-spencer}
    E(\gamma) = \Tr \, \Tilde{F}_{\Tilde{A}}\wedge \Tilde{F}_{\Tilde{A}}- \Tr \, F_{A}\wedge F_{A},
\end{equation}
provided that right-hand side is perpendicular to $\ker E$. To verify this, we first note that integration by parts shows that the kernel is given by
\begin{equation*}
    \ker E = \{ \varphi \in \Lambda^{2,2}(X) : d \varphi = 0, \quad \bar{\partial}^\dagger \partial^\dagger \varphi = 0 \}.
\end{equation*}
From here it is easy to see that $i \partial \bar{\partial} \chi \perp \ker E$, and so we may solve
\begin{equation} \label{kodaira-spencer2}
    E(\gamma) = i \partial \bar{\partial} \chi
\end{equation}
which gives a solution to \eqref{kodaira-spencer}. Next, a well-known integration by parts argument \cite{FLY, PicardTopics} shows that equation \eqref{kodaira-spencer2} implies $d \gamma = 0$. Therefore
\begin{equation*}
\partial \bar{\partial} \bar{\partial}^\dagger \partial^\dagger \gamma = \Tr \, \Tilde{F}_{\Tilde{A}}\wedge \Tilde{F}_{\Tilde{A}}- \Tr \, F_{A}\wedge F_{A},
\end{equation*}
and we let
\begin{equation*}
    R_2[\tilde{A},A] := -i \bar{\partial}^\dagger \partial^\dagger \gamma.
\end{equation*}
In other words
\begin{equation} \label{defn-R2}
    R_2[\tilde{A},A] = -i \bar{\partial}^\dagger \partial^\dagger E^{-1} (\Tr \, \Tilde{F}_{\Tilde{A}}\wedge \Tilde{F}_{\Tilde{A}}- \Tr \, F_{A}\wedge F_{A})
\end{equation}
We take the real part if necessary so that $R_2[\tilde{A},A]$ is a real $(1,1)$-form. 


\subsection{Deformation of gauge connection}
For both Bott-Chern and Aeppli deformations of the metric, we choose the same deformation path for the gauge connection from the parameter space $H^{1}(X,\End \, E)$, such that the deformed curvature $\Tilde{F}_{\Tilde{A}}^{0,2}$ automatically satisfies $\Tilde{F}_{\Tilde{A}}^{0,2}=0$. Recall that assumption (BG4) ensures that $H^1(X,\End \, E)$ locally parametrizes solutions to the equation $\{ F^{0,2} = 0 \}$.
\smallskip
\par Let $E$ be a vector bundle over $X$ with reference metric $h$. Let $\bar{\D}_\alpha$ be a smooth family of holomorphic structures on $E$ varying with a parameter $\alpha$. This means $\bar{\D}_\alpha: \Gamma(E) \rightarrow \Omega^{0,1}(E)$ satisfies $\bar{\D}_\alpha^2=0$. By the Chern correspondence, the operator $\bar{\D}_\alpha$ uniquely determines a connection $D_\alpha = \D_\alpha + \bar{\D}_\alpha$ which is metric compatible with respect to $h$.
\smallskip
\par We now change the metric by the relation
\[
  \tilde{h} = h e^u
\]
where $u \in \Gamma({\rm End} \, E)$ satisfies $u^\dagger = u$, and $\dagger$ denotes the adjoint with respect to $h$. The adjoint of an endomorphism $\gamma \in \Gamma({\rm End} \, E)$ satisfies $\langle \gamma x, y \rangle_h = \langle x, \gamma^\dagger y \rangle_h$, i.e. $\gamma^\dagger = (h \gamma h^{-1})^*$. We remark that according to our conventions,
\[
\langle x,y \rangle_{\tilde{h}} = \langle e^u x, y \rangle_h.
\]
The Chern correspondence associates to the holomorphic structure $\bar{\D}_\alpha$ a unique $\tilde{h}$-metric compatible connection
\[
D_{\alpha,u} = \D_{\alpha,u} + \bar{\D}_\alpha
\]
with
\begin{equation} \label{eq:GaugePath}
\D_{\alpha,u} = \D_\alpha + e^{-u} \D_\alpha e^u.
\end{equation}
Equation \eqref{eq:GaugePath} may be derived by imposing metric compatibility \eqref{metric-comp} with respect to $\tilde{h}$. At a fixed parameter $\alpha$, we may open-up a holomorphic frame for $\bar{\D}_\alpha$ in which we see
\[
\D_{\alpha,u} =\partial+  \tilde{h}^{-1} \partial \tilde{h}, \quad \bar{\D}_\alpha = \bar{\partial}
\]
and
\[
\tilde{F}_{\alpha,u} = \bar{\partial} (\tilde{h}^{-1} \partial \tilde{h}), \quad \tilde{h} = h e^u.
\]
For the setup with additional spurious gauge field on the tangent bundle, its deformation follows similarly. Namely, from two parameters $\alpha_1$ and $\alpha_2$ smoothly parametrizing holomorphic structures $\bar{\D}_{\alpha_1}$ on $E$ and $\bar{\D}_{\alpha_2}$ on $T^{1,0}X$, a pair of reference metrics $(h,g)$ on $(E,T^{1,0}X)$ gives metric compatible connections $(E,D_{\alpha_1})$ and $(T^{1,0}X,D_{\alpha_2})$. The ansatz is then
\[
\D_{\alpha_1,u_1} = \D_{\alpha_1} + e^{-u_1} \D_{\alpha_1} e^{u_1}, \quad \D_{\alpha_2,u_2} = \D_{\alpha_2} + e^{-u_2} \D_{\alpha_2} e^{u_2},
\]
for $u_1 \in \Gamma({\rm End} \, E)$ and $u_2 \in \Gamma({\rm End} \, T^{1,0}X)$. For our setup below, the reference metric $h$ is Hermitian-Yang-Mills \eqref{ref-metric} on a holomorphic stable bundle $E \rightarrow (X,\omega_{\rm CY})$. On $T^{1,0}X$, the reference metric $g$ is the Calabi-Yau metric.

\subsection{Deformation of complex structure}
In this paper, we fix the complex structure along both deformation paths. The difficulty comes from the fact that the deformation of complex structure changes $\bar{\partial}$ which couples to both deformation of metric and gauge connection. Due to this, one usually needs more effort to construct a natural deformation ansatz. Therefore, we leave it to future investigations.

\section{The Setup for the Implicit Function Theorem}
Once we have determined our ansatz of the deformation paths, we seek zeros of the remaining equations. Therefore, we separate coordinates as 
\begin{equation}
    \mathrm{X}=(\alpha',[\text{cohomology class}]),\quad \mathrm{Y}=(\text{parameter within equivalence class}) \nonumber
\end{equation}
and construct the map $\F$ such that
\begin{equation}
    \F(\mathrm{X},\mathrm{Y})=\qty[\text{remaining heterotic equations}] \nonumber
\end{equation}
and look for the zeros of the map $\F(\textrm{X},\textrm{Y})$. Then, if the zero locus exists, we will have constructed heterotic solutions, and $\textrm{X}$ will be the desired coordinate of the parameter space. 

Note that since we start the deformation path from a K\"ahler solution, $\F(0,0)=0$. We then compute the linearization at $\alpha'=0$, $\eval{D_{\textrm{Y}}\F}_{0}$ and show that it is an isomorphism. Then by the implicit function theorem, there exists $\epsilon>0$, such that for any $\alpha'<\epsilon$, there exists solutions $\textrm{Y}(\textrm{X})$ such that $\F(\textrm{X},\textrm{Y}(\textrm{X}))=0$. In other words, for each $\alpha'$, there is a neighborhood of the starting K\"ahler solution parameterized by $\textrm{X}$, which produces solutions $(\Tilde{\omega},\Tilde{A})$ to the heterotic system. 

This strategy was used previously by Li and Yau \cite{LY05} and Andreas and Garcia-Fernandez \cite{AG2012} to construct solutions to the Strominger system. Our work brings in the various cohomology classes (Bott-Chern and Aeppli) into the picture, and the explicit ansatzes \eqref{eq:BCPath}, \eqref{eq:APath} show how these cohomology classes can be used as parameters.

The rest of the paper proceeds to establish the validity of such application of the implicit function theorem on our two deformation paths corresponding to the Bott-Chern case and the Aeppli case.

\subsection{Bott-Chern case} 
The deformation ansatz reads by \eqref{eq:BCPath} and \eqref{eq:GaugePath},
    \begin{align*} |\Omega|_{\Tilde{\omega}} \, \Tilde{\omega}^2&=|\Omega|_{\omega} \, \omega^2+\mathfrak{b}+\Theta,\\
       D_{\alpha_1,u_1} &= \D_{\alpha_1,u_1} + \bar{\D}_{\alpha_1}, \quad \D_{\alpha_1,u_1} = \D_{\alpha_1} + e^{-u_1} \D_{\alpha_1} e^{u_1}.
    \end{align*}
The parameters $\mathfrak{b}$ and $\alpha_1$ should represent cohomology classes. Recall that part of our hypotheses is that $\bar{\D}_{\alpha_1}$ is a smooth family of holomorphic structures on $E$ parametrized by $\alpha_1 \in H^{1}({\rm End} \, E)$ \cite{Huybrechts}. To set this up as a Banach space, we use the Hodge theorem to identify cohomology with the space of harmonic representatives. Let
\begin{equation*}
    \mathbb{H}^{p,p} = \{ \Theta \in \Lambda^{p,p}(X,\mathbb{R}) : \Delta_{\omega_{\rm CY}} \Theta = 0 \}.
\end{equation*}
By Hodge theory this is a finite-dimensional vector space. By the $\partial \bar{\partial}$-lemma, the $(2,2)$-form $\Theta$ being $d$-exact is equivalent to $\partial \bar{\partial}$-exact. We will use deformation coordinates $\textrm{X}$ and $\textrm{Y}$ given by
    \begin{align}
        \textrm{X}&=(\alpha',\mathfrak{b},\alpha_{1})\in \mathbb{R}\times \mathbb{H}^{2,2} \times H^{1}({\rm End} \, E), \nonumber\\
        \textrm{Y}&=(\Theta,u_1)\in {\rm Im} \, d \cap \Lambda^{2,2}(X,\mathbb{R})\times \Gamma( {\rm End}_0  E), \nonumber
    \end{align}
    where
    \[
     \Gamma({\rm End}_0  E) = \{ u \in \Gamma({\rm End} \,E) : u^\dagger = u, \ \ {\rm Tr} \, u = 0 \}.
    \]
    The ${\rm End}_0$ is needed since shifts $u \mapsto u+C \, {\rm id}_E$ by constant multiples of the identity are not seen by this ansatz. As this ansatz is automatically conformally balanced, the map $\F_{\text{BC}}$ is constructed by including the remaining equations in the system.
    \begin{align}
        &\F_{\text{BC}}(\textrm{X},\textrm{Y})=\mqty[i\partial\bar{\partial}\Tilde{\omega}-\alpha'(\Tr \, \Tilde{R}_{\Tilde{\omega}}\wedge \Tilde{R}_{\Tilde{\omega}} - \Tr \, \Tilde{F}_{\alpha_1,u_1}\wedge \Tilde{F}_{\alpha_1,u_1})\\| \Omega|_{\Tilde{\omega}} e^{u/2} \, \Tilde{\omega}^2\wedge i\Tilde{F}_{\alpha_1,u_1} \, e^{-u/2}].
    \end{align}
   Solutions to the system are then cut out by $\F_{\text{BC}}(\textrm{X},\textrm{Y})=0$. Our initial data is a K\"ahler Ricci-flat metric $\omega=\omega_{\rm CY}$ and Donaldson-Uhlenbeck-Yau connection $A$, and so $\F_{\text{BC}}(0,0)=0$.
\smallskip   
\par To apply the implicit function theorem on Banach spaces, we take $k \geq 2$, $0<\gamma<1$, and setup the domain and codomain to be
   \begin{align}
           \F_{\text{BC}}&: \mathbb{R}\times \mathbb{H}^{2,2} \times H^1({\rm End} \, E) \times C^{k+2,\gamma}({\rm Im} \, d \cap \Lambda^{2,2}) \times C^{k+2,\gamma}({\rm End}_0 E) \nonumber \\
            &\qquad\qquad\qquad\qquad\qquad\qquad\qquad\longrightarrow C^{k,\gamma}({\rm Im} \, d \cap \Lambda^{2,2}) \times C^{k,\gamma}(V(E)). \nonumber
           \end{align}
    
    It is worth noting that we have further constrained the image of $\F_{\text{BC}}$ to the following subspace
    \begin{equation} \label{V-defn}
        V(E)=\{v\in\Lambda^{6}(X,\End E):\int_{X}\Tr \, v=0, \ v^{\dagger}=v\}.
    \end{equation}
   This will be needed when showing that the linearized operator is surjective. We also note that the spaces involved are all Banach spaces. For example, $C^{k,\gamma}({\rm Im} \, d \cap \Lambda^{2,2})$ is a Banach space (e.g. \cite{marshall}).

It is easy to check that the image of $\F_{\text{BC}}$ is contained in $\partial \bar{\partial} \Lambda^{1,1} \times V(E)$. The first row is contained in the image of $\partial \bar{\partial}$, since we assume $c_2(E)=c_2(X)$ and we may apply the $\partial \bar{\partial}$-lemma. For containment of the second row in $V(E)$, it suffices to show that the following two conditions hold:
\begin{enumerate}
    \item Vanishing integration of trace:
    \begin{equation*}
        \int_{X}\Tr (|\Omega|_{\Tilde{\omega}} \,\Tilde{\omega}^{2}\wedge i\Tilde{F})=0.
    \end{equation*}
This uses that $c_1(E) = 0$ so that ${\rm Tr} \, \tilde{F} = d \gamma$, and the ansatz is such that $|\Omega|_{\Tilde{\omega}}\, \Tilde{\omega}^{2}$ is closed.
    \item Self-adjoint:
    \begin{equation*}
        (|\Omega|_{\Tilde{\omega}}e^{u/2} \, \Tilde{\omega}^{2}\wedge i\Tilde{F} e^{-u/2} )^{\dagger}=|\Omega|_{\Tilde{\omega}} e^{u/2} \, \Tilde{\omega}^{2}\wedge i\Tilde{F} \, e^{-u/2}
    \end{equation*}
    This property is the reason for conjugating by $e^{u/2}$ in the definition of $\F_{\text{BC}}$. Note that here the adjoint $\dagger$ is taken with respect to the reference metric $h$ while $\tilde{F}$ is the Chern curvature of the deformed metric given by $\tilde{h} = h e^u$ with $(e^u)^\dagger = e^u$. To show
    \[
      \langle e^{u/2} i \tilde{F} e^{-u/2}x , y \rangle_h = \langle x , e^{u/2} i \tilde{F} e^{-u/2} y \rangle_h,
    \]
   one can use $(i \tilde{F})^{\dagger_{\tilde{h}}} = i \tilde{F}$ and $\langle x,y \rangle_{\tilde{h}} = \langle e^u x,y \rangle_h$.
\end{enumerate}

    To establish the validity of the implicit function theorem and conclude the existence of the heterotic solution in neighborhood of $(\textrm{X},\textrm{Y})=(0,0)$, we need to check the invertibility of $D_{\textrm{Y}}\F_{\text{BC}}$ at $(\textrm{X},\textrm{Y})=(0,0)$,
    \begin{align}
        \eval{D_{\textrm{Y}}\F_{\text{BC}}}_{(0,0)}&: C^{k+2,\gamma}({\rm Im} \, d \cap \Lambda^{2,2}) \times C^{k+2,\gamma}({\rm End}_0 E)  \longrightarrow C^{k,\gamma}({\rm Im} \, d \cap \Lambda^{2,2}) \times C^{k,\gamma}(V(E)) \nonumber\\
        &\eval{D_{\textrm{Y}}\F_{\text{BC}}}_{(0,0)}(\dot{\Theta},\dot{u}_1)=\mqty[L^{\text{BC}}_1 & 0\\ C & L_2]\mqty[\dot{\Theta}\\\dot{u}_1], \label{eq:DF-BC}
    \end{align}
    and prove that $\eval{D_{\textrm{Y}}\F_{\text{BC}}}_{(0,0)}$ is an isomorphism. 

    \subsection{Aeppli case}
    Recall the setup for the reference geometry: $\omega$ is a K\"ahler Ricci-flat metric with Chern connection $\theta$, $A$ is a Hermitian-Yang-Mills connection over $E \rightarrow (X,\omega)$ and $\beta \in \Omega^{1,1}(X,\mathbb{R})$ is an $i \partial \bar{\partial}$ potential for ${\rm Tr} \, R_\theta^2 - {\rm Tr} \, F_A^2$. Let $\alpha_1 \in H^1({\rm End} \, E)$ parametrize a smooth family of holomorphic structures on $E$ and $\alpha_2 \in H^1({\rm End} \, T^{1,0}X)$ parametrize a smooth family of holomorphic structures on $T^{1,0}X$.

    The deformation ansatz reads by \eqref{eq:APath} and \eqref{eq:GaugePath}
    \begin{align*}
      \Tilde{\omega}&=\omega+\mathfrak{a}+b(\xi)+\alpha'\qty(-R_{2}[\Tilde{A},A]+R_{2}[\Tilde{\theta},\theta] + \beta), \\
       D_{\alpha_i,u_i} &= \D_{\alpha_i,u_i} + \bar{\D}_{\alpha_i}, \quad \D_{\alpha_i,u_i} = \D_{\alpha_i} + e^{-u_i} \D_{\alpha_i} e^{u_i},
    \end{align*}
    where $\tilde{A}$ and $\tilde{\theta}$ are given by $D_{\alpha_1,u_1} = d + \tilde{A}$ and $D_{\alpha_2,u_2} = d + \tilde{\theta}$ and
    \begin{equation*}
b(\xi) = (1+J) d \xi \in \Omega^{1,1}(X,\mathbb{R}), \quad \xi \in \Lambda^{1}(X,\mathbb{R}).
      \end{equation*}
    The deformation parameters $\textrm{X}$ and $\textrm{Y}$ read,
    \begin{align}
        \textrm{X}&=(\alpha',\mathfrak{a}, \alpha_1 , \alpha_2) \in \mathbb{R}\times \mathbb{H}^{1,1} \times H^1(\End \, E)\times H^1(\End \, T^{1,0}X),\nonumber \\
        \textrm{Y}&=(\xi,u_1,u_2)\in ({\rm Im} \, d^\dagger) \cap \Lambda^{1}(X,\mathbb{R}) \times \Gamma({\rm End}_0  E) \times \Gamma({\rm End}_0  T^{1,0}X). \nonumber
    \end{align}

    \begin{rmk}
We note that the anomaly cancellation relation \eqref{anomaly-eqn3} is automatically implemented by the ansatz, which can also be found in \cite{MGF22}. 
      \end{rmk}
    
      \begin{rmk}
        The ansatz for the $(1,1)$-form $b$ ensures that it is in the image of $\partial \oplus \bar{\partial}$. Therefore
 \begin{equation*}
     \underline{a}(\tilde{\omega},\tilde{A},\tilde{\theta}) = [\omega]_A + [\mathfrak{a}]_A
 \end{equation*}
 which is the main property of this ansatz. We write $\mathfrak{a} \in \mathbb{H}^{1,1}$ using harmonic representatives of de Rham cohomology, but we could equivalently write $\mathfrak{a} \in H^{1,1}_A$ as we are on a K\"ahler manifold and Aeppli and de Rham cohomology coincide (Lemma 5.15 in \cite{DGMS1975}).
\end{rmk}
    
    The map $\F_{\text{A}}$ is taken to be
    \begin{align}
        \F_{\text{A}}(\textrm{X},\textrm{Y})=\mqty[\star_{\omega} d(|\Omega|_{\Tilde{\omega}} \Tilde{\omega}^2)\\|\Omega|_{\Tilde{\omega}} e^{u_1/2}\Tilde{\omega}^2\wedge i\Tilde{F} e^{-u_1/2} - d_1 \,  \tilde{\omega}^3 \otimes {\rm id} \\|\Omega|_{\Tilde{\omega}} e^{u_2/2} \Tilde{\omega}^2\wedge i\Tilde{R} e^{-u_2/2} - d_2  \, \tilde{\omega}^3 \otimes {\rm id}],
    \end{align}
    where
    \begin{align}
        d_1 &= \frac{1}{{\rm rk} \, E} \frac{\int_X |\Omega|_{\tilde{\omega}} \tilde{\omega}^2 \wedge {\rm Tr} \, i \tilde{F}}{\int_X |\Omega|_{\tilde{\omega}} \tilde{\omega}^3} \nonumber\\
 d_2 &= \frac{1}{3} \frac{\int_X |\Omega|_{\tilde{\omega}} \tilde{\omega}^2 \wedge {\rm Tr} \, i \tilde{R}}{\int_X |\Omega|_{\tilde{\omega}} \tilde{\omega}^3}. \nonumber
    \end{align}
    The domain and codomain are
    \begin{equation*}
\F_{\text{A}}: \mathcal{X} \times \mathcal{Y} \rightarrow \mathcal{Z},
    \end{equation*}
    with
    \begin{align*}
      \mathcal{X} &=\mathbb{R}\times  \mathbb{H}^{1,1} \times H^1(\End \, E)\times H^1(\End \, T^{1,0}X) \\
      \mathcal{Y} &= C^{k+2,\gamma}[({\rm Im} \, d^\dagger) \cap \Lambda^{1}(X,\mathbb{R})] \times C^{k+2,\gamma}({\rm End}_0  E) \times C^{k+2,\gamma}({\rm End}_0  T^{1,0}X) \\
      \mathcal{Z} &=C^{k,\gamma}[({\rm Im} \, d^\dagger) \cap \Lambda^{1}(X,\mathbb{R})] \times C^{k,\gamma}(V(E)) \times C^{k,\gamma}(V(T^{1,0}X)).
      \end{align*}
    Here $d^\dagger_\omega = - \star_\omega d \star_\omega$ is with respect to the background Calabi-Yau metric, and $C^{k,\gamma}[({\rm Im} \, d^\dagger) \cap \Lambda^{r}(X,\mathbb{R})]$ is well-known to be a Banach space \cite{marshall}. The inclusion of the $d_1$, $d_2$ terms is so that the image of the bottom rows of $\F$ are contained in $V(E)$ and $V(T^{1,0}X)$ as before \eqref{V-defn}, meaning that the integral of the trace is equal to zero.

The set $\F_{\text{A}}(\textrm{X},\textrm{Y}) = 0$ is our solution space. This is because if
\begin{equation*}
    d (|\Omega|_{\tilde{\omega}} \tilde{\omega}^2)=0, \quad |\Omega|_{\Tilde{\omega}}\Tilde{\omega}^2\wedge i\Tilde{F} = d_1 \,  \tilde{\omega}^3 \otimes {\rm id},
\end{equation*}
   and then upon taking the trace and integrating we see that
   \begin{equation*}
\int_X |\Omega|_{\Tilde{\omega}}\Tilde{\omega}^2\wedge i {\rm Tr} \, \Tilde{F} =0
   \end{equation*}
   since $c_1(E)=0$ implies ${\rm Tr} \, \Tilde{F} = d \gamma$. Therefore $d_1=0$ and similarly $d_2=0$ on shell.
\smallskip
 \par We should comment on why $\F_{\textrm{A}}$ is a differentiable map. This is a standard argument as everything is explicit, except for differences
\begin{equation*}
 R_2[\tilde{A}+h,A] - R_2[\tilde{A},A].   
\end{equation*}
First, we note that for connections $A,\tilde{A} \in C^{k+1,\gamma}$ and $h \in C^{k+1,\gamma}(\Lambda^1({\rm End} \, E))$, then $R[\tilde{A}+h,A] \in C^{k+2,\gamma}$. This is from the expression \eqref{defn-R2} and elliptic regularity
\begin{equation*}
    E^{-1}: C^{k,\gamma} \rightarrow C^{k+4,\gamma}
\end{equation*}
for the 4th order elliptic operator $E$. Next, we set $L_{\tilde{A}} h = {d \over dt} \big|_{t=0} ({\rm Tr} \, F_{\tilde{A}+th} \wedge F_{\tilde{A}+th})$, and explain why $R_2[\tilde{A},A]$ is differentiable in the unknown $\tilde{A}$, namely
\begin{equation} \label{differentiability}
\lim_{h \rightarrow 0} {\| R_2[\tilde{A}+h,A] - R_2[\tilde{A},A] - \bar{\partial}^\dagger \partial^\dagger E^{-1}(L_{\tilde{A}} h) \|_{C^{k+2,\gamma}} \over \| h \|_{C^{k+1,\gamma}}} = 0.
\end{equation}
For this, we use \eqref{defn-R2} to obtain
\begin{align}
& \ \| R_2[\tilde{A}+h,A] - R_2[\tilde{A},A] - \bar{\partial}^\dagger \partial^\dagger E^{-1}(L_{\tilde{A}} h) \|_{C^{k+2,\gamma}} \nonumber\\
&= \| \bar{\partial}^\dagger \partial^\dagger E^{-1} ({\rm Tr} \, F_{\tilde{A}+h} \wedge F_{\tilde{A}+h} - {\rm Tr} \, F_{\tilde{A}} \wedge F_{\tilde{A}} - L_{\tilde{A}} h ) \|_{C^{k+2,\gamma}}.    \nonumber  
\end{align}
To obtain bounds on $E^{-1}$, we use the standard elliptic estimate
\begin{equation*}
    \| \gamma \|_{C^{4+k,\gamma}} \leq C \| E(\gamma) \|_{C^{k,\gamma}}, \quad \gamma \in (\ker E)^\perp.
\end{equation*}
Thus
\begin{align}
& \ \| R_2[\tilde{A}+h,A] - R_2[\tilde{A},A] - \bar{\partial}^\dagger \partial^\dagger E^{-1}(L_{\tilde{A}} h) \|_{C^{k+2,\gamma}} \nonumber\\
&\leq C \| ({\rm Tr} \, F_{\tilde{A}+h} \wedge F_{\tilde{A}+h} - {\rm Tr} \, F_{\tilde{A}} \wedge F_{\tilde{A}} - L_{\tilde{A}} h) \|_{C^{k,\gamma}}    \nonumber\\
&\leq C \bigg\| \int_0^1 {d \over dt} ({\rm Tr} \, F_{\tilde{A} + th} \wedge F_{\tilde{A} + th}) dt - {d \over dt}\bigg|_{t=0} ({\rm Tr} \, F_{\tilde{A} + th} \wedge F_{\tilde{A} + th}) \bigg\|_{C^{k,\gamma}}. \nonumber
\end{align}
Since
\begin{equation*}
    {d \over dt} {\rm Tr} \, (F_t)^2 = 2 {\rm Tr} \, F_t d_{A_t} h, \quad A_t = \tilde{A} + th,
\end{equation*}
we have
\begin{align}
& \ \| R_2[\tilde{A}+h,A] - R_2[\tilde{A},A] - \bar{\partial}^\dagger \partial^\dagger E^{-1}(L_{\tilde{A}} h) \|_{C^{k+2,\gamma}} \nonumber\\
&\leq C \bigg\| \int_0^1 ({\rm Tr} F_t d_{A_t} - {\rm Tr} F_{\tilde{A}} d_{\tilde{A}}) h \, dt \bigg\|_{C^{k,\gamma}} \nonumber\\
&= C \bigg\| \int_0^1 \int_0^1 {d \over ds} ({\rm Tr} \, F_s d_{A_s}) h ds \, dt \bigg\|_{C^{k,\gamma}} \nonumber
\end{align}
where $A_s = \tilde{A} + tsh$. It follows from here that 
\begin{equation*}
     \| R_2[\tilde{A}+h,A] - R_2[\tilde{A},A] - \bar{\partial}^\dagger \partial^\dagger E^{-1}(L_{\tilde{A}} h) \|_{C^{k+2,\gamma}} \leq C \| h \|_{C^{k+1,\gamma}}^2,
\end{equation*}
where the constant $C$ depends on $\| \tilde{A} \|_{C^{k+1,\gamma}}$. The limit \eqref{differentiability} follows, and so does differentiability of $\tilde{\omega}(\mathfrak{a},b,\tilde{A},\tilde{\theta})$ in the argument $\tilde{A}$ (and similarly for $\tilde{\theta}$). 
\smallskip
\par    With this setup in place, $\textrm{X}$ is the moduli coordinate and we will compute the linearization of $\F$ with respect to $\textrm{Y}$, denoted $D_{\textrm{Y}}\F_{\text{A}}$, at $\alpha'=0$:
    \begin{align}
        \eval{D_{\textrm{Y}}\F_{\text{A}}}_{0}&: C^{k+2,\gamma}[{\rm Im} \, d^\dagger \cap \Lambda^{1}(X,\mathbb{R})] \times C^{k+2,\gamma}({\rm End}_0  E) \times C^{k+2,\gamma}({\rm End}_0  T^{1,0}X) \nonumber\\
        &\qquad\longrightarrow C^{k,\gamma}[{\rm Im} \, d^\dagger \cap \Lambda^{1}(X,\mathbb{R})]  \times C^{k,\gamma}(V(E)) \times C^{k,\gamma}(V(T^{1,0}X)) \nonumber\\
        &\eval{D_{\textrm{Y}}\F_{\text{A}}}_{0}(\dot{\xi},\dot{u}_1,\dot{u}_2)=\mqty[L^{\text{A}}_1 & 0 & 0\\ C_{1} & L_2 & 0 \\ C_2 & 0 & L_3]\mqty[\dot{\xi}\\\dot{u}_1\\\dot{u}_2]. \label{eq:DF-A}
    \end{align}
    We will prove the invertibility of this operator and conclude the existence of heterotic solutions $\F_{\text{A}}(\text{X},\text{Y}(\text{X}))=0$ near
    \begin{equation*}
    \F_{\text{A}}(0,0)=0    
    \end{equation*}
    by the implicit function theorem. The point $(0,0)$ is in the zero locus of $\F_{\text{A}}$ since our initial data is a K\"ahler Ricci-flat metric $\omega=\omega_{\rm CY}$, a Donaldson-Uhlenbeck-Yau connection $A$, and the Chern connection of the K\"ahler Ricci-flat metric on $T^{1,0}X$.

\section{Calculations of $D\F$}
In this section, we proceed to show the invertibility of $D_{\textrm{Y}}\F$ for both cases and establish the validity of the application of the implicit function theorem to conclude the existence of the heterotic solutions near a K\"ahler solution. 

\subsection{Calculation of $D_{\textrm{Y}}\F_{\text{BC}}$}
Consider a deformation path $(\Theta(s),u_1(t))$ along $X=0$ with $(\Theta(0),u_{1}(0))=(0,0)$ since the starting point is a K\"ahler solution. Then
\begin{align*}
    |\Omega|_{\tilde{\omega}(s)} \tilde{\omega}^2(s) &= |\Omega|_\omega \omega^2 + \Theta(s) \nonumber\\
    D_{u(t)} &= \D_0 + e^{-u(t)} \D_0 e^{u(t)} + \bar{\D}_0 
\end{align*}
and we will compute the partial derivatives $\partial_s|_0$, $\partial_t|_0$ of $\F$ along this path.

\subsubsection{Invertibility of $\eval{D_{\textrm{Y}}\F_{\text{BC}}}_{(0,0)}$ at $\alpha'=0$}
For the implicit function theorem to hold, we need $D_{\textrm{Y}}\F_{\text{BC}}$ to be an isomorphism at $(\textrm{X},\textrm{Y})=(0,0)$, hence $D_{\textrm{Y}}\F_{\text{BC}}$ should have an inverse $(D_{\textrm{Y}}\F_{\text{BC}})^{-1}$. From \eqref{eq:DF-BC}, it must be
\begin{equation}\label{eq:bc-f-1}
    (D_{\textrm{Y}}\F_{\text{BC}})^{-1}=\mqty[(L^{\text{BC}}_{1})^{-1} & 0 \\ -(L^{\text{BC}}_{1})^{-1}CL_2^{-1} & L_2^{-1}].
\end{equation}
Therefore, we only need to establish the invertibility of $L^{\text{BC}}_{1}$ and $L_2$ respectively. 

\subsubsection{$L^{\text{BC}}_{1}$}
To calculate $L^{\text{BC}}_{1}$, we need to understand the variation of the hermitian metric $\Tilde{\omega}$ along the Bott-Chern deformation path \eqref{eq:BCPath}. This is provided by the following lemma.

Notation: In this sub-subsection, we write $\delta_{s=0}$ as $\delta$ and $\delta\Theta=\dot\Theta$ for simplicity. 

\begin{lemma}[Variation of $\Tilde{\omega}$ under Bott-Chern deformation]
On the Bott-Chern deformation path \eqref{eq:BCPath} given by $|\Omega|_{\Tilde{\omega}}\Tilde{\omega}^2=|\Omega|_{\omega} \omega^2+\Theta$,  then
\begin{equation}\label{eq:w-var-bc}
    \delta\Tilde{\omega}=\frac{1}{2|\Omega|_{\omega}}\Lambda_{\omega}\delta\Theta.
\end{equation}
\end{lemma}

This formula is well-known and can be found in e.g. \cite{CPY2022, PPZ1}. This calculation can be completed by various methods. We provide here a proof which uses the Lefschetz operator commutator relations.

\begin{proof}
    First, we find the variation of $|\Omega|_{\Tilde{\omega}}$ at $\alpha'=0$, using the local expression \eqref{eq:localOmega}, 
\begin{align*}
    \delta| \Omega|^2_{\Tilde{\omega}}&=\Omega\bar{\Omega}\delta\qty(\det \Tilde{g}_{\alpha\bar{\beta}})^{-1}\nonumber \\
    2|\Omega|_{\omega}\delta|\Omega|_{\Tilde{\omega}}&=-\Omega\bar{\Omega}\det(g_{\alpha\bar{\beta}})^{-2}\delta\qty(\det \Tilde{g}_{\alpha\bar{\beta}})\nonumber \\
    &=-\Omega\bar{\Omega}\det(g_{\alpha\bar{\beta}})^{-1}g^{\bar{\beta}\alpha}\delta \Tilde{g}_{\alpha\bar{\beta}}\nonumber \\
    &=-|\Omega|^2_{\omega}g^{\bar{\beta}\alpha}\delta \Tilde{g}_{\alpha\bar{\beta}}\nonumber \\
    &=-|\Omega|^2_{\omega}(\Lambda_{\omega}\delta\Tilde{\omega}). \nonumber \\
    \implies \delta|\Omega|_{\Tilde{\omega}}&=-\frac{|\Omega|_{\omega}}{2}(\Lambda_{\omega}\delta\Tilde{\omega}).
\end{align*}
Then the variation of $(|\Omega|_{\Tilde{\omega}}\Tilde{\omega}^2)$ follows,
\begin{align}
    \delta(| \Omega|_{\Tilde{\omega}}\Tilde{\omega}^2)&=\qty(\delta(|\Omega|_{\Tilde{\omega}})\omega^2+|\Omega|_{\omega}\delta(\Tilde{\omega}^2)) \nonumber\\
    \delta\Theta&=|\Omega|_{\omega}\qty(-\tfrac{1}{2}\qty(\Lambda_{\omega}\delta\Tilde{\omega})\omega^2+2\omega\wedge\delta\Tilde{\omega})\label{eq:variation}
\end{align}

We can rewrite \eqref{eq:variation} using $L_{\omega}$ and contract both side with $\omega$, and then invoke the Lefschetz identities \eqref{eq:KahlerId} to obtain,
\begin{align*}
     \delta\Theta&=|\Omega|_{\omega}\qty(-\tfrac{1}{2}\qty(\Lambda_{\omega}\delta\Tilde{\omega})L_{\omega}\omega+2L_{\omega}\delta\Tilde{\omega}) \\
    \implies\Lambda_{\omega}\delta\Theta&=|\Omega|_{\omega}\qty(-\tfrac{1}{2}\qty(\Lambda_{\omega}\delta\Tilde{\omega})\Lambda_{\omega}L_{\omega}\omega+2\Lambda_{\omega}L_{\omega}\delta\Tilde{\omega}) \nonumber\\
    &=|\Omega|_{\omega}\qty(-\tfrac{1}{2}\qty(\Lambda_{\omega}\delta\Tilde{\omega})L_{\omega}\Lambda_{\omega}\omega+2L_{\omega}\Lambda_{\omega}\delta\Tilde{\omega}) \nonumber \\
    &\qquad-|\Omega|_{\omega}\qty(-\tfrac{1}{2}\qty(\Lambda_{\omega}\delta\Tilde{\omega})H\omega+2H\delta\Tilde{\omega}) \nonumber\\
    &=|\Omega|_{\omega}\qty(\tfrac{1}{2}(-2n+2+4)(\Lambda_{\omega}\delta\Tilde{\omega})\omega+2(n-2)\delta\Tilde{\omega}) \nonumber \\
   &=2|\Omega|_{\omega}\delta\Tilde{\omega}.
\end{align*}
By rearranging, we arrive at the desired variation of $\Tilde{\omega}$ as in \eqref{eq:w-var-bc}.
\end{proof}

Then $L^{\text{BC}}_{1}$ follows from \eqref{eq:w-var-bc} immediately:
\begin{align*}
    L^{\text{BC}}_{1}\dot{\Theta}(0)&=i\partial\bar{\partial}\delta\Tilde{\omega}=\frac{i\partial\bar{\partial}}{2|\Omega|_{\omega}}\qty(\Lambda_{\omega}\delta\Theta)=-\frac{\partial\partial^{\dagger}_{\omega}\dot{\Theta}}{2|\Omega|_{\omega}}=-\frac{\Delta_{\omega}\dot{\Theta}}{2|\Omega|_{\omega}}, \nonumber\\
    &\implies L^{\text{BC}}_{1}=-\frac{1}{2|\Omega|_{\omega}}\Delta_{\omega}.
\end{align*}
Here we have taken advantage of the K\"ahler background at $s=0$ to use the K\"ahler identity $[\Lambda_{\omega},\bar{\partial}]=-i\partial^{\dagger}_{\omega}$ and write $\Delta_{\omega}=\partial\partial^{\dagger}_{\omega}+\partial^{\dagger}_{\omega}\partial$ as $\Theta$ and $\dot{\Theta}$ are $\partial\bar{\partial}$-exact. 

We proceed to establish the invertibility of $L^{\text{BC}}_{1}$, which follows from Hodge theory. Indeed, $\Delta_\omega = d d^\dagger + d^\dagger d$ by the K\"ahler condition, $|\Omega|_\omega$ is constant, and our operator is
\begin{align*}
    L_1^{\text{BC}} &: {\rm Im} \, d \cap \Lambda^{2,2} \rightarrow {\rm Im} \, d \cap \Lambda^{2,2} \nonumber\\
    L_1^{\text{BC}} &= -{1 \over 2 |\Omega|_\omega} \Delta_\omega
\end{align*}
which is invertible by Hodge theory. Indeed, $L_1^{\text{BC}}$ is injective since if $\Delta \theta =0$ and $\theta = d \alpha$ then integration by parts shows $\theta =0$. Furthermore, $L_1^{\text{BC}}$ is surjective since given $\eta \in \Lambda^{2,2}$ with $\eta = d \beta$ one can solve $\Delta \theta = \eta$ with $\theta = d \alpha$, and $\theta \in \Lambda^{2,2}$ follows since $\Delta$ preserves type on K\"ahler manifolds. Elliptic regularity (e.g. \cite{marshall}) provides an inverse in H\"older spaces via the Schauder estimate
\begin{equation*}
\| \Theta \|_{C^{k+2,\gamma}} \leq C \| \Delta_\omega \Theta \|_{C^{k,\gamma}}, \quad \Theta \in {\rm Im} \, d \cap \Lambda^{2,2}.
\end{equation*}
Thus
\begin{equation*}
    (L_1^{\text{BC}})^{-1}: C^{k,\gamma}({\rm Im} \, d \cap \Lambda^{2,2}) \rightarrow C^{k+2,\gamma}({\rm Im} \, d \cap \Lambda^{2,2})
\end{equation*}
and $L^{\text{BC}}_{1}$ is invertible.

\subsubsection{$L_2$}
To calculate $L_2$, we need to understand the deformed curvature $\Tilde{F}$ on the deformation path \eqref{eq:GaugePath} with holomorphic bundle parameter set to $\alpha_1=0$. Written in a holomorphic frame for the initial holomorphic bundle $(E,\bar{\D})$, the curvature is
\[
\Tilde{F} = \bar{\partial} (\tilde{h}^{-1} \partial \tilde{h}), \quad \tilde{h} = h e^u.
\]
We parametrize $u(t)$ with $u(0)=0$ and write $\delta_{t=0}$ as $\delta$ and $\delta u=\dot{u}$. We note
\[
\tilde{h}^{-1} \partial \tilde{h} = e^{-u} h^{-1} \partial (h e^u) = h^{-1} \partial h + e^{-u} \D e^u,
\]
where $\D e^u = \partial e^u + [A,e^u]$ and $A = h^{-1} \partial h$. From this, we can derive the variation of $\tilde{F}$:
\begin{equation*}
    \delta\Tilde{F} = \bar{\D} \D \dot{u}.
\end{equation*}
Using this and $i\Lambda_\omega F =0$ at $u=0$, the linearization of the HYM equation becomes
\begin{align*}
    L_{2}\dot{u}&=\delta(|\Omega|_{\omega}e^{u/2} \omega^2\wedge i\Tilde{F} e^{-u/2}) \nonumber \\
    &=\delta\qty( e^{u/2} i\Lambda_{\omega}\Tilde{F} e^{-u/2}\otimes|\Omega|_{\omega}\tfrac{\omega^3}{3!}) \nonumber \\
    &=i\Lambda_{\omega}\delta\Tilde{F} \otimes|\Omega|_{\omega}\tfrac{\omega^3}{3!} \nonumber \\
    &=i\Lambda_{\omega}\qty(\bar{\D}\D \dot{u})\otimes|\Omega|_{\omega}\tfrac{\omega^3}{3!}\nonumber \\
    \implies L_{2}&=i\Lambda_{\omega}\qty(\bar{\D}\D)\otimes|\Omega|_{\omega}\tfrac{\omega^3}{3!}.
\end{align*}
We obtain
\begin{align*}
L_2 &: C^{k+2,\gamma}({\rm End}_0 E) \rightarrow C^{k,\gamma}(V(E)), \nonumber\\
L_2 (\dot{u}) &= - g^{\mu \bar{\nu}} \bar{\D}_{\bar{\nu}} \D_\mu \dot{u} \otimes |\Omega|_{\omega}\tfrac{\omega^3}{3!}.
\end{align*}
We now proceed to investigate the kernel of $L_2$. Take $\gamma \in \Gamma({\rm End}_0 E)$ with $\gamma\in\Ker(\Lambda_{\omega}\bar{\D}\D)$.
Integrating by parts and using the balanced condition, we can deduce from 
\begin{equation} \label{L2-kernel}
\langle \Lambda_{\omega}\bar{\D}\D \gamma,\gamma \rangle=0
\end{equation} 
that $\gamma$ is holomorphic, $\bar{\D}\gamma=0$. See Appendix \ref{app:kernellambda} for details.

The assumption in our setup is that $E \rightarrow (X,\omega)$ is stable. Then the stability implies that $E$ is simple, hence any holomorphic endomorphism section $\gamma$ of $\End E$ is proportional to the identity, i.e. $\gamma = \lambda\cdot\id_{E}$. Together with the constraint $\Tr \, \gamma=0$, our kernel of $L_2$ must be zero for $\lambda=0$, as expected. 

Next, we show that $L_2$ is surjective. For this, we use the Fredholm alternative $\Im(L_2)=(\Ker L_2^{\dagger})^{\perp}$ (e.g. Appendix of \cite{DonaldsonKronheimer}) where the adjoint $\dagger$ is with respect to the $L^2$ inner product \eqref{End-inner-prod}. If $\gamma \in \Ker L_2^{\dagger}$ then \eqref{L2-kernel} holds. It follows as above that $\gamma = \lambda \cdot \id_E$. Since $\gamma \in V(E)$, then $\lambda=0$ and $\gamma=0$. Therefore $\Ker L_2^{\dagger} = \{ 0\}$ and $L_2$ is surjective. Elliptic regularity gives the invertibility of $L_2$.
\begin{equation*}
    L_2^{-1} : C^{k,\gamma}(V(E)) \rightarrow C^{k+2,\gamma}({\rm End}_0 E).
\end{equation*}

Altogether, by \eqref{eq:bc-f-1} the invertibility of $L^{\text{BC}}_1$ and $L_2$ is sufficient for the invertibility of $D_Y\F$, and hence the linearization $D_{\textrm{Y}}\F_{\text{BC}}$ is indeed an isomorphism. Then the application of the implicit function theorem for Banach spaces concludes Theorem \ref{thm-bc}.

\subsection{Calculation of $D_{\textrm{Y}}\F_{\text{A}}$}
We now return to the Aeppli deformation setup. Recall that here the $\textrm{Y}$ deformations consist of
\[
\xi \in {\rm Im} \, d^\dagger \cap \Lambda^1, \quad u_1 \in \Gamma({\rm End}_0 E), \quad u_2 \in \Gamma({\rm End}_0 T^{1,0}X).
  \]
The $\textrm{X}$ deformations include $\alpha' \in \mathbb{R}$ and the moduli coordinates $\mathbb{H}^{1,1}$, $H^1({\rm End} \, E)$ and $H^1({\rm End} \, T^{1,0}X)$.

\subsubsection{Invertibility of $\eval{D_{\textrm{Y}}\F_{\text{A}}}_{0}$ at $\alpha'=0$}
Similarly to the Bott-Chern case, we need $D_{\textrm{Y}}\F_{\text{A}}$ to be an isomorphism at $(\textrm{X},\textrm{Y})=(0,0)$ and show $D_{\textrm{Y}}\F$ has an inverse $(D_{\textrm{Y}}\F_{\text{A}})^{-1}$. From \eqref{eq:DF-A}, we have,
\begin{equation*}
(D_{\textrm{Y}}\F_{\text{A}})^{-1}=\mqty[(L^{\text{A}}_1)^{-1} & 0 & 0 \\ -(L^{\text{A}}_1)^{-1}C_{1}L_{2}^{-1} & L_2^{-1} & 0 \\
    -(L^{\text{A}}_1)^{-1}C_{2}L_{3}^{-1} & 0 & L_3^{-1} ].
\end{equation*}
Therefore, we only need to establish the invertibility of $L^{\text{A}}_{1}$, $L_2$, and $L_3$ respectively. To calculate linearizations, we use a deformation path $(\xi(s),u_1(t_1),u_2(t_2)) \in \textrm{Y}$ along $\textrm{X}=0$ with starting point the K\"ahler solution $(\xi(0),u_{1}(0),u_2(0))=(0,0,0)$.

\subsubsection{$L^{\text{A}}_{1}$}
To calculate $L^{\text{A}}_{1}$, we need to understand the deformation of metric on the Aeppli deformation path \eqref{eq:APath}. In particular, we can write the local expression of $\Tilde{\omega}$ as,
\begin{equation*}
\tilde{g}_{\alpha\bar{\beta}}=g_{\alpha\bar{\beta}}-i\mathfrak{a}_{\alpha\bar{\beta}}-ib_{\alpha\bar{\beta}}-i\mathcal{O}(\alpha'),
\end{equation*}
where we denote $\mathfrak{a}=\mathfrak{a}_{\alpha\bar{\beta}}dx^{\alpha\bar{\beta}}$ and $b=b_{\alpha\bar{\beta}}dx^{\alpha\bar{\beta}}$. We parameterize $b(s)$ and write $\delta_{s=0}$ as $\delta$ and $\delta\Tilde{\omega}=\delta b=\dot{b}$ for simplicity. By the ansatz for $b$, we have $\dot{b} = (1+J) d \dot{\xi}$. By previous calculation \eqref{eq:variation}, we have at $s=0$
\begin{align}
    \delta\qty(\star_\omega d\qty(|\Omega|_{\Tilde{\omega}}\Tilde{\omega}^2))&=\star_{\omega}d\qty(\delta\qty(|\Omega|_{\Tilde{\omega}}\Tilde{\omega}^2))\nonumber \\
    &=\star_{\omega}d\qty[|\Omega|_{\omega}\qty(-\tfrac{1}{2}\qty(\Lambda_{\omega}\delta\Tilde{\omega})\omega^2+2\omega\wedge\delta\Tilde{\omega})]\nonumber .
\end{align}
Using that $d \omega=0$ and $|\Omega|_\omega$ is constant, we may rearrange this as
\begin{equation*}
   \delta\qty(\star_\omega d\qty(|\Omega|_{\Tilde{\omega}}\Tilde{\omega}^2)) =  2 |\Omega|_{\omega} \star_{\omega} (\omega \wedge d T \dot{b} ), \quad \delta \tilde{\omega} = \dot{b}
  \end{equation*}
  where
  \begin{equation*}
T\dot{b} = -\tfrac{1}{4}\qty(\Lambda_{\omega}\dot{b})\omega+\dot{b}.
\end{equation*}
Thus
\begin{equation*}
L^{\text{A}}_{1}\dot{\xi} = 2 |\Omega|_{\omega} \star_{\omega} \qty( \omega \wedge d T (1+J) d \dot{\xi} ).
\end{equation*}
We now invoke the identity of \cite{WuPL}, which states that on a compact K\"ahler threefold there holds
\begin{equation*}
\star_{\omega} \qty( \omega \wedge d T (1+J) d \dot{\xi} ) = \Delta_d \dot{\xi}, \quad \dot{\xi} \in {\rm Im} \, d^\dagger \cap \Lambda^1.
\end{equation*}
Here $\Delta_d = d d^\dagger + d^\dagger d$. We refer to \cite{WuPL} for the full calculation. Therefore
\begin{align*}
    L_1^{\text{A}} &: {\rm Im} \, d \cap \Lambda^{1} \rightarrow {\rm Im} \, d \cap \Lambda^{1} \nonumber\\
    L_1^{\text{A}}  &= 2 |\Omega|_{\omega} \Delta_d .
\end{align*}
The constant factor $2 |\Omega|_\omega$ is irrelevant, and standard Hodge theory of the Laplacian on forms lets us conclude that $L^{\text{A}}_{1}$ admits an inverse
\begin{equation*}
(L^{\text{A}}_{1})^{-1}:  C^{k,\gamma}({\rm Im} \, d^\dagger \cap \Lambda^{1}) \rightarrow C^{k+2,\gamma} ({\rm Im} \, d^\dagger \cap \Lambda^{1})
\end{equation*}
as needed.

\subsubsection{$L_2$ and $L_3$}
The $L_2$ and $L_3$ are the same compared to the Bott-Chern case, hence we can conclude their invertibility respectively.

Therefore, the implicit function theorem for Banach spaces can be applied to conclude the existence of solutions of heterotic systems near the small $\epsilon$-neighborhood of the reference K\"ahler solution, i.e. Theorem \ref{thm-a}.

\section{Further discussions}
We have established the existence of local coordinates near a K\"ahler solution with two different deformation paths. Several questions arise naturally: 
\begin{enumerate}
\item The duality between $H^{2,2}_{\text{BC}}(X,\mathbb{R})$ and $H^{1,1}_{\text{A}}(X,\mathbb{R})$ implies that their dimensions are equal and suggests that the constructed local coordinates can be related by some coordinate transformation in the tangent space of the moduli space. However, the current work introduces a spurious gauge field $\theta$ in the Aeppli case. After this paper was posted, an analogous result was obtained in the Aeppli case \cite{WuPL} without introducing $\theta$, which allows the Bott-Chern and Aeppli deformations to be considered on the same set of equations.

    \item We have established the local coordinates of the heterotic moduli near a K\"ahler solution by taking advantage of many nice properties of K\"ahler manifold. Can we establish the existence of local coordinates on a general heterotic solution? It seems natural to extend our setup to an arbitrary solution of the heterotic system (1). However, our setup is limited as we have invoked some nice properties of K\"ahler threefolds. In particular, the linearized operators $L_1^{\rm BC}$ and $L_1^{\rm A}$ are proportional to $\Delta_\omega$ and $(d+d^\dagger)$ in the K\"ahler case. In the non-K\"ahler case,  the operators $L_1^{\rm BC}$ and $L_1^{\rm A}$ are more complicated and it is unclear how to show that these are invertible. We have also used the $\partial\bar{\partial}$-lemma at various points in the argument (such as using de Rham cohomology to represent harmonic forms in Bott-Chern or Aeppli cohomology), and we refer to e.g. \cite{AT2017} for various characterizations of the $\partial \bar{\partial}$-lemma involving Bott-Chern and Aeppli cohomology.

    
    \item Throughout our paper, we have fixed the complex structure on the complex manifold $X$. Therefore, technically speaking, the constructed local coordinates are not true local coordinates of the moduli space but provide local coordinates to a subspace of the moduli instead. In particular, it is the local coordinates of the subspace of the heterotic moduli where the complex structure is invariant. Hence, a further extension to the full moduli is needed for future investigations. We remark that in general, the existence of Strominger solutions on compact non-K\"ahler manifolds is not stable by deformation of the complex structure; see \cite{IU19}.
    
    \item While our setup involves varying the metric and gauge field $(g,A)$, other investigations of the moduli e.g. \cite{CdlOMc,MGF17} introduce an additional field called the $B$-field and vary $(g,B,A)$. One motivation for adding the $B$-field is to give the moduli space the structure of a complex manifold \cite{CdlOMcS,MGF-JDG}. It would be interesting to further investigate $B$-field and obtain the complexified description. 
\end{enumerate}

\newpage
\appendix
\section{The physical system}\label{app:justforchern}

The equations of heterotic supergravity are derived by an expansion in a small parameter denoted $\alpha'$. In this model, the fields are a metric tensor $g$, the curvature $F$ of a connection on a vector bundle $E$, a scalar function $\Phi$, and a 3-form $H$. The 3-form is encoded into the geometry by producing a connection $\nabla^{\rm H}$ which has Christoffel symbols
\begin{equation}\label{hull-connection}
\Gamma^{\rm H}{}_i{}^k{}_j = \Gamma^{\rm LC}{}_i{}^k{}_j + {1 \over 2} H_i{}^k{}_j.
\end{equation}
The action \cite{BdR} to first order in $\alpha'$ is given by
\begin{equation*}
S = \int e^{-2 \Phi} \bigg[ R + 4 |\nabla \Phi|^2 - {1 \over 2} |H|^2 + \alpha' ({\rm Tr} \, |F|^2 - {\rm Tr} \, |R^{\rm H}|^2) \bigg] d {\rm vol}_g + O(\alpha'^2).
  \end{equation*}
The principle of least action states that the physical equations should correspond to stationary points of the action. 
\smallskip
\par Special critical points of the action involve extra structure on the manifold. Certain non-K\"ahler structures in complex geometry were identified by Strominger \cite{Strominger} as satisfying the equations of supersymmetry. In this setup, $X$ is a complex manifold with holomorphic volume form $\Omega$ and hermitian metric $\omega= i g_{\mu \bar{\nu}} dx^{\mu \bar{\nu}}$, and the associated fields are given by
  \begin{equation*}
H = i (\partial-\bar{\partial})\omega + O(\alpha'^2), \quad \Phi = - {1 \over 2} \log |\Omega|_\omega + O(\alpha'^2).
  \end{equation*}
  The supersymmetry constraints are
  \begin{equation} \label{susy}
d(|\Omega|_\omega \omega^2) = O(\alpha'^2), \quad F^{0,2}=O(\alpha'), \quad F \wedge \omega^2 = O(\alpha')
    \end{equation}
    and the heterotic Bianchi identity is
  \begin{equation} \label{hull-anomaly}
i \partial \bar{\partial} \omega = \alpha' ( {\rm Tr} \, R^{\rm H} \wedge R^{\rm H} - {\rm Tr} \, F \wedge F ) + O(\alpha'^2).
\end{equation}
The consistency of these equations with conformal invariance of the sigma model and anomaly cancellation was established by Hull \cite{H1986}. It is well-known in the string theory literature (see \cite{CGGL} and more recent work e.g. \cite{AMP23,OS2014b}) that solutions to \eqref{susy}, \eqref{hull-anomaly} are critical points of the action $S$.
 \smallskip
 \par As an aside, we note that with $H=i(\partial-\bar{\partial})\omega$ and the conventions considered here, then \eqref{hull-connection} defines the Hull connection (which does not preserve the complex structure) while 
\begin{equation}\label{bismut-connection}
    \hat{\Gamma}_i{}^k{}_j = \Gamma^{\rm LC}{}_i{}^k{}_j - {1 \over 2} H_i{}^k{}_j.
\end{equation}
is the Strominger-Bismut connection which does preserve the complex structure.
  \smallskip
\par From the perspective of string theory, these constraints should not be taken as literal equality, as they are not valid to higher order in $\alpha'$. That said, as a mathematical problem we may look for solutions where the right-hand side of the supersymmetry constraints \eqref{susy} is exactly zero. Such non-K\"ahler geometries are interesting as generalizations of Calabi-Yau metrics to non-K\"ahler Calabi-Yau threefolds, and these equations have applications to the geometrization of conifold transitions \cite{CPY24, CGPY, FLY}. However, the constraint \eqref{hull-anomaly} is problematic as a literal equality without the higher order $O(\alpha'^2)$ contributions, since $ {\rm Tr} \, R^{\rm H} \wedge R^{\rm H}$ is generally not a $(2,2)$ form.
\smallskip
\par The solutions constructed in this paper solve the mathematical equations \eqref{system-of-eqn} near a K\"ahler background. The equations \eqref{system-of-eqn} are almost the same as the truncated physical system, except that ${\rm Tr} \, R \wedge R$ is computed with respect to the Chern connection instead of the Hull connection $\nabla^{\rm H}$. We now explain why in the setting of this paper, this is not a problem as these solutions do satisfy the physical system to this order in $\alpha'$. This is because applying the implicit function theorem over a K\"ahler background gives a smooth path $g_{\alpha'}$ in the parameter $\alpha'$ with $g|_{\alpha'=0} = g_{\rm CY}$. In other words,
  \begin{equation} \label{g-orderbyorder}
    g_{\alpha'} = g_{\rm CY} + O(\alpha').
  \end{equation}
  It follows that these solutions satisfy the estimate $H=O(\alpha')$ and so
  \begin{equation*}
 R^{\rm H} - R^{\rm Ch} = \mathcal{O}(H, \nabla H) = O(\alpha').
\end{equation*}
Thus \eqref{hull-anomaly} is satisfied and solutions to the heterotic Bianchi identity with Chern connection are not distinguishable from solutions with Hull connection at this order in $\alpha'$. 
\smallskip
\par We also considered another setup in this paper where the heterotic Bianchi identity is taken with extra connection $\theta$ on $T^{1,0}X$.
\begin{equation} \label{4form-BI}
i \partial \bar{\partial} \omega = \alpha' ({\rm Tr} \, R_\theta \wedge R_\theta - {\rm Tr} \, F \wedge F ).
\end{equation}
Here we applied the implicit function theorem over a K\"ahler background to obtain a path $\theta_{\alpha'}$ such that $\theta_{\alpha'} \rightarrow \Gamma^{\rm LC}_{g_{\rm CY}}$ as $\alpha' \rightarrow 0$. It follows that
\begin{equation*}
R_{\theta} - R^{\rm H}_{g_{\rm CY}} = O(\alpha').
\end{equation*}
Combining this with the fact that the metric $g$ is Calabi-Yau at zeroth order \eqref{g-orderbyorder}, it follows that \eqref{hull-anomaly} is satisfied. 
\smallskip
\par In summary, the mathematical equations solved by Li and Yau \cite{LY05} and Andreas and Garcia-Fernandez \cite{AG2012} (see \cite{CPY2022} for follow-up work) by inverse function theorem near a K\"ahler Calabi-Yau solution solve the physical system. The present work deforms these solutions by the implicit function theorem, and so both setups considered in this paper satisfy the physical system. From the perspective of string theory, exact solutions to the truncated equations are not needed as solving the physical equations by large radius expansion about a K\"ahler Calabi-Yau is well-known in the string theory literature by works of Witten-Witten \cite{WittenWitten} and Melnikov-Minasian-Sethi \cite{MMS}. 
\smallskip
\par We make one further comment on the difference between the setup of the physical system in string theory and the setup considered in the current paper. A fundamental equation in our approach is the 4-form Bianchi identity \eqref{4form-BI}, while in string theory the anomaly cancellation relation is typically written as a 3-form equation
\begin{equation*}
    H = dB + \alpha' (CS[A]-CS[\theta])
\end{equation*}
where $CS[A]$, $CS[\theta]$ are Chern-Simons local 3-forms and $B$ is a local 2-form which transforms on overlaps in such a way that $H$ is a well-defined global 3-form. In string theory, variations of the 3-form field strength $H$ are related to variations of the 2-form $B$ field \cite{AGS,CdlOMc}. The approach to the moduli problem in the present work is to vary the 3-form $H$ directly by variations of the metric $\omega$ via $H = i(\partial-\bar{\partial}) \omega$.

\section{Kernel of $\Lambda_{\omega}\bar{\D}\D$}\label{app:kernellambda}
Let $(E,\bar{\D}) \rightarrow X$ be a holomorphic vector bundle, and let $(\omega,\Omega)$ be a conformally balanced structure so that $d(|\Omega|_\omega \omega^2)=0$. We define the inner product on sections of $\End E$ as
\begin{equation} \label{End-inner-prod}
    \langle \gamma_1,\gamma_2 \rangle =\int_{X}\Tr(\gamma_1\gamma_2^{\dagger})|\Omega|_{\omega}\frac{\omega^3}{3!}, \quad \forall \gamma_1, \gamma_2\in\Gamma(\End E),
\end{equation}
where the adjoint $\dagger$ is with respect to a metric $h$. Recall that $D = \D + \bar{\D}$ is a metric compatible connection with respect to $h$, namely the Chern connection of $h$ on the holomorphic bundle $(E,\bar{\D})$. We will show that if $\gamma \in \Gamma({\rm End} \, E)$ with $\gamma^\dagger = \gamma$ solves $\gamma\in\Ker(\Lambda_{\omega}\bar{\D}\D)$, then $\gamma$ is holomorphic. Using the HYM equation $i\Lambda_{\omega}F_{A}=0$, we have
\begin{align*}
    0&= \langle i\Lambda_{\omega} \bar{\D}\D \gamma,\gamma \rangle \nonumber\\
    &= - \int_{X}\Tr(i\Lambda_{\omega}( \D_{\mu}\bar{\D}_{\bar{\nu}}\gamma-[F_{\mu\bar{\nu}},\gamma])dx^{\mu\bar{\nu}})\cdot\gamma^{\dagger}|\Omega|_{\omega}\frac{\omega^3}{3!}\nonumber\\
    &= - \int_{X}i\Lambda_{\omega}\Tr (\D_{\mu}\bar{\D}_{\bar{\nu}}\gamma\cdot \gamma^{\dagger})dx^{\mu\bar{\nu}}|\Omega|_{\omega}\frac{\omega^3}{3!}=: - I_1.
\end{align*}
For the evaluation of $I_1$, integration by part yields
\begin{align*}
    I_1&=\int_{X}i\Lambda_{\omega}\Tr(\D_{\mu}\bar{\D}_{\bar{\nu}}(\gamma\gamma^{\dagger}))dx^{\mu\bar{\nu}}|\Omega|_{\omega}\frac{\omega^3}{3!}\nonumber\\
    &\ -\int_{X}i\Lambda_{\omega}\Tr(\bar{\D}_{\bar{\nu}}\gamma\cdot\D_{\mu}\gamma^{\dagger})dx^{\mu\bar{\nu}}|\Omega|_{\omega}\frac{\omega^3}{3!}\nonumber\\
    &\ -\int_{X}i\Lambda_{\omega}\Tr(\D_{\mu}\gamma\cdot\bar{\D}_{\bar{\nu}}\gamma^{\dagger})dx^{\mu\bar{\nu}}|\Omega|_{\omega}\frac{\omega^3}{3!}\nonumber\\
    &\ -\int_{X}i\Lambda_{\omega}\Tr(\gamma\cdot\D_{\mu}\bar{\D}_{\bar{\nu}}(\gamma^{\dagger}))dx^{\mu\bar{\nu}}|\Omega|_{\omega}\frac{\omega^3}{3!}.
\end{align*}
Since $\gamma$ is self-adjoint, the above yields
\begin{equation*}
\langle \D \gamma, \D \gamma \rangle  = -I_1 + \frac{1}{2}\int_{X} i\Lambda_{\omega}\Tr(\D_{\mu}\bar{\D}_{\bar{\nu}}(\gamma\gamma^{\dagger}))dx^{\mu\bar{\nu}}|\Omega|_{\omega}\frac{\omega^3}{3!}.
\end{equation*}
Setting $I_1=0$ and commuting $\D$ with the trace, we obtain
\begin{align*}
    \langle \bar{\D}\gamma,\bar{\D}\gamma \rangle 
    &=\frac{1}{2}\int_{X}i\Lambda_{\omega}\partial\bar{\partial}\Tr(\gamma\gamma^{\dagger})|\Omega|_{\omega}\frac{\omega^3}{3!}\nonumber\\
    &=\frac{3}{2}\int_{X}i\partial\bar{\partial}\Tr(\gamma\gamma^{\dagger}) \wedge |\Omega|_{\omega}\frac{\omega^2}{3!}=:\frac{3}{2\cdot 3!}I_2,
\end{align*}
where the last equality holds as for any function $f$, we have the identity
\begin{equation*}
(\Lambda_{\omega} i \partial\bar{\partial}f) \omega^{n}=n i \partial\bar{\partial}f\wedge\omega^{n-1}.
\end{equation*}
Integrating by parts again, we have
\begin{equation*}
    I_2=\int_{X}i d ( \bar{\partial}\Tr(\gamma\gamma^{\dagger}) \wedge |\Omega|_{\omega}\omega^2)\nonumber\\
    + \int_{X}i \bar{\partial}\Tr(\gamma\gamma^{\dagger}) \wedge d( |\Omega|_{\omega}\omega^2)
    =0,
\end{equation*}
where the first term vanishes by Stoke's theorem and compactness of base manifold $X$, and the other vanishes due to the balanced condition $d(|\Omega|_{\omega}\omega^2)=0$. Therefore, $\bar{\D} \gamma=0$, or in other words, $\gamma$ is holomorphic.

\end{document}